\documentclass{amsart}

\usepackage{amssymb}
\usepackage{multirow}
\newtheorem{theorem}{Theorem}[section]
\newtheorem{lemma}[theorem]{Lemma}

\theoremstyle{definition}
\newtheorem{definition}[theorem]{Definition}

\theoremstyle{remark}
\newtheorem{remark}[theorem]{Remark}

\newcommand\relphantom[1]{\mathrel{\phantom{#1}}}
\numberwithin{equation}{section}

\numberwithin{equation}{section}

\begin{document}

\title[Galerkin Approximations for Stochastic Fractional Wave Equations]{Galerkin Finite Element Approximations for Stochastic Space-Time  Fractional Wave Equations}

\author{Yajing Li}
\address{School of Mathematics and Statistics, Gansu Key Laboratory of Applied Mathematics and Complex Systems, Lanzhou University, Lanzhou 730000, P.R. China }
\email{liyajing11@lzu.edu.cn}
\author{Yejuan Wang}
\address{School of Mathematics and Statistics, Gansu Key Laboratory of Applied Mathematics and Complex Systems, Lanzhou University, Lanzhou 730000, P.R. China}
\email{wangyj@lzu.edu.cn}

\author{Weihua Deng}
\address{School of Mathematics and Statistics, Gansu Key Laboratory of Applied Mathematics and Complex Systems, Lanzhou University, Lanzhou 730000, P.R. China}
\email{dengwh@lzu.edu.cn}



\begin{abstract}
The traditional wave equation models wave propagation in an ideal conducting medium. For characterizing the wave propagation in inhomogeneous media with frequency dependent power-law attenuation, the space-time fractional wave equation appears; further incorporating the additive white Gaussian noise coming from many natural sources leads to the stochastic space-time fractional wave equation. This paper discusses the Galerkin finite element approximations for the stochastic space-time fractional wave equation forced by an additive space-time white noise. We firstly  discretize the space-time additive noise, which introduces a modeling error and results in a regularized stochastic space-time fractional wave equation; then the regularity of the regularized equation is analyzed. For the discretization in space, the finite element approximation is used and the definition of the discrete fractional Laplacian is introduced. We derive the mean-squared $L^2$-norm priori estimates for the modeling error and for the approximation error to the solution of the regularized problem; and the numerical experiments are performed to confirm the estimates. 
For the time-stepping, we calculate the analytically obtained Mittag-Leffler type function.

\end{abstract}

\maketitle
\pagestyle{myheadings}
\thispagestyle{plain}
\markboth{Y. J. LI, Y. J. WANG, W. H. DENG}{GALERKIN APPROXIMATIONS FOR STOCHASTIC FRACTIONAL WAVE EQUATIONS}


\section{Introduction}
In an ideal conducting medium, the propagation of wave is well governed by the classical wave equation. But, most of the time, because of the complex  inhomogeneity of the conducting medium, the wave propagation with frequency-dependent attenuation has been observed in a wide range of areas, including acoustics, viscous dampers in seismic isolation of buildings, structural vibration, seismic wave propagation \cite{Chen04, Meerschaert, Szabo94}. More often, the attenuated wave propagation exhibits a power-law relation between attenuation and frequency, e.g., acoustical wave propagation in lossy media \cite{Szabo94}. The exponent of power-law generally ranges from $0$ to $2$, which belongs to anomalous attenuation \cite{Chen04}, and the competing model for describing the corresponding attenuated wave propagation is the space-time fractional wave equations, in which the fractional derivatives play the role of characterizing the power-law behavior.

Over the last few decades, one of the main reasons that the fractional calculus attracts wide interests is because it can effectively characterize the ubiquitous power-law phenomena. Even in the kingdom of diffusion, the central pillars of anomalous diffusion, including sub-diffusion and super-diffusion, long-range dependence (LRD), $1/f$ noise, and L\'{e}vy statistics, have intrinsic power-law structures \cite{Eliazar}; and the corresponding fractional diffusion equations have been much better investigated than the power law wave equation \cite{Bouchaud,Deng,Jin3,Jin,Solomon,Zeng}. For the diffusion-wave equations, there are also some discussions, including analysis, algorithm, and applications; the fundamental solutions and their properties are considered in \cite{Agarwal,Mainardi,Mainardi3,Sakamoto}; and different kinds of numerical methods and approximation schemes have also been developed, e.g., the finite difference method \cite{Sun,Zhang}, spectral method \cite{Bhrawy}, and finite element method \cite{Esen,Jin4}, etc.

For the practically physical system, the different stochastic perturbations are coming from many natural sources; sometimes, they can not be ignored and we need to incorporate them to the corresponding deterministic model; and then the stochastic differential equations (SDEs) are produced.  In more recent decades, the initial and/or boundary value problems for SDEs have been extensively studied theoretically \cite{Chen1, Hu, Peter, Kulinich} and numerically \cite{Armstrong, BABUSKA, du, Kovacs, Kovacs1, Milstein, Yan}. However, it seems that there are less literatures related to the theoretical analysis or numerical approximation of stochastic wave equations with fractional derivative or driven by fractional noise. In \cite{Caithamer}, the stochastic wave equation forced  by a class of fractional noises is explicitly solved and the upper and lower bounds on both the large and small deviations of several sup norms associated with the solution are given. The 1-d stochastic wave equation driven by a fractional Brownian sheet is discussed in \cite{Sardanyons}, and a Young integration theory is developed.
This paper is considering the stochastic space-time fractional wave equation driven by an additive space-time white noise; to the best of our knowledge, so far there are no research results on the effective numerical schemes of this kind of equation. For filling this gap, we provide the finite element approximation for the following initial boundary value problem with $1<\alpha<2$
and $\frac{1}{2}<\beta\leq 1$:
\begin{equation} \label{eq0}
\left\{ \begin{array}
 {l@{\quad} l}
\displaystyle \partial_t^{\alpha}u(t,x)+(-\Delta)^\beta u(t,x)=\frac{\partial^2W(t,x)}{\partial t \partial x},
 \quad 0<t<T,~0<x<1,\\
 \\
 u(t,0)=u(t,1)=0, \quad  0<t<T,\\
 \\
  u(0,x)=v_1(x),~ \partial_tu(0,x)=v_2(x),\quad  0<x<1,
 \end{array}
 \right.
\end{equation}
where $\partial_t^{\alpha}$ denotes the left-sided Caputo fractional derivative of order $\alpha$
with respect to $t$; $(-\Delta)^\beta$ is the fractional Laplacian, the definition of which is based on the spectral
decomposition of the Dirichlet Laplacian, as adopted in \cite{Nochetto}; and $W(t,x)$ represents an
infinite dimensional Brownian motion. It is worth noting that
if $\beta=1$, then $(-\Delta)^\beta$ will reduce to the Laplace operator, which has been studied in \cite{Jin4}
for the corresponding deterministic case.
Moreover, when $\alpha=\beta=1$ in \eqref{eq0}, error estimates of solutions have been obtained in \cite{du}.
 Obviously, this is an interesting and important generalization of Galerkin finite element methods for the stochastic
parabolic equations with a classical derivative in \cite{du}.
As we know, one of the main difficulties in using the Galerkin finite element method to solve the variational formulation of \eqref{eq0} is how to discretize the fractional Laplacian. For this, motivated by the definition of the discrete Laplace operator in \cite{Jin}, we firstly propose the definition of the discrete fractional Laplacian in this paper.

The paper is organized as follows. In Section \ref{sec:2}, we first introduce some basic definitions, notations
and necessary preliminaries, then prove the stability of the equation and discretize the space-time additive noise, which introduces a modeling error and results in a regularized stochastic space-time fractional wave equation; and the convergence orders of the modeling error are well established, which has a transition point $\alpha=\frac{3}{2}$, i.e., the convergence order is $\alpha-\frac{1}{2}-\epsilon$ for $\alpha \in (1,\frac{3}{2}]$ and 1 for $\alpha \in (\frac{3}{2},2]$, where $\alpha$ is the order of the time fractional derivative and $\epsilon$ can be any sufficiently small positive constant. In Section \ref{sec:3}, the regularity of the regularized stochastic space-time fractional wave equation is discussed; the finite element scheme is presented and the corresponding error estimates are detailedly derived. In Section \ref{sec:4}, the numerical experiments are performed to confirm the convergence orders of the modeling error and finite element approximations to the regularized equation. We conclude the paper with some remarks in the last section.

\section{Preliminaries}\label{sec:2}
In this section, we collect useful facts on the Mittag-Leffler function, present a representation of the solution of problem \eqref{eq0}, give the stability results of the solution and consider the approximation of noise.

\subsection{Mittag-Leffler function.}
We shall frequently use the Mittag-Leffler function $E_{\alpha,\beta}(z)$ defined as follows:
\[E_{\alpha,\beta}(z)=\sum\limits_{k=0}^{\infty}\frac{z^k}{\Gamma(k\alpha+\beta)},\quad z\in \mathbb{C}.\]
The Mittag-Leffler function $E_{\alpha,\beta}(z)$ is a two-parameter family of entire functions in $z$ of order $\alpha^{-1}$ and type 1 \cite[p. 42]{Kilbas}. It generalizes the exponential function in the sense that $E_{1,1}(z)=e^z$. For later use, we collect some results in the next lemma; see \cite{Kilbas, Podlubny}.
\begin{lemma}\label{le2.1}
Let $0<\alpha<2$ and $\beta\in\mathbb{R}$ be arbitrary. We suppose that $\mu$ is an arbitrary real number
 such that $\frac{\pi\alpha}{2}<\mu< \min(\pi,\pi\alpha)$.
Then there exists a constant $C=C(\alpha,\beta,\mu)>0$ such that
\begin{equation}\label{eq2.1}
|E_{\alpha,\beta}(z)|\leq \frac{C}{1+|z|}, \quad\quad \mu\leq|\arg(z)|\leq \pi.
\end{equation}
Moreover, for $\lambda>0$, $\alpha>0$ and positive integer $m \in \mathbb{N}$, we have
\begin{equation}\label{eq2.2}
\frac{d^m}{d t{^m}} E_{\alpha,1}(-\lambda^\beta t^\alpha)=
-\lambda^\beta t^{\alpha-m}E_{\alpha,\alpha -m+1}(-\lambda^\beta t^\alpha), \quad t > 0,
\end{equation}
and
\begin{equation}\label{eq2.3}
\frac{d}{d t}\left(tE_{\alpha,2}(-\lambda^\beta t^\alpha)\right)=
E_{\alpha,1}(-\lambda^\beta t^\alpha), \quad t \geq 0.
\end{equation}
\end{lemma}
\subsection{Solution representation.}
To discuss the stability of solutions and the approximation of noise, we first need some notations, and some basic definitions for the Caputo fractional derivative and fractional Laplacian. Let $H=L^2(0,1)$ with $\|\cdot\|$ to denote the norm in $H$ and
$(\cdot, \cdot)$ to denote the inner product of ${H}$. It is well known that the Laplacian $-\Delta=-\frac{\partial^2}{\partial x^2}$ has eigenpairs $(\lambda_k,e_k)$ with $\lambda_k=k^2\pi^2$, $e_k=\sqrt{2}\sin k\pi x$, $k=1,2,3,\ldots$ subject to the homogeneous Dirichlet boundary conditions on $(0,1)$, i.e., $e_k(0)=e_k(1)=0,$ and $-\Delta e_k=\lambda_k e_k, k=1,2,3,\ldots.$ The set $\{e_k\}_{k=1}^\infty$ forms an orthonormal basis in $H$.
For $q\geq 0$, we denote by $\dot{H}^q\subset H$ the Hilbert space induced by the norm
\[|u|_q^2=\sum\limits_{k=1}^{\infty}\lambda_k^q(u,e_k)^2.\]
Then, for any $u\in \dot{H}^{2\beta},$  we have
\[(-\Delta)^\beta u=\sum\limits_{k=1}^{\infty}\lambda_k^\beta(u,e_k)e_k,\]
which is the fractional Laplacian we consider in this work.

The definitions of the fractional integral and left-sided Caputo fractional derivative of order $\alpha$ with respect to $t$ can be found in \cite[p. 91]{Kilbas} or \cite[p. 78]{Podlubny}.
\begin{definition}\label{def1}
The fractional integral of order $\alpha> 0$ with the lower limit $0$ for a function $u$ is defined as
\[I^\alpha_t u(t) =\frac{1}{\Gamma(\alpha)}\int^t_{0}(t-s)^{\alpha-1}u(s)ds, \quad t > 0, \]
where $\Gamma(\cdot)$ is the gamma function.
\end{definition}
\begin{definition}\label{def2}
The Caputo fractional derivative of order $\alpha > 0$ with the lower limit $0$ for a function $u$ is defined as
\[\partial^\alpha_t u(t) =\frac{1}{\Gamma(n-\alpha)}\int^t_{0}
(t-s)^{n-\alpha-1}\frac{\partial ^nu(s)}{\partial s^n}ds, \quad t > 0, \quad 0\leq n-1 < \alpha < n,\]
where the function $u(t)$ has absolutely continuous derivatives up to order $n-1$.
\end{definition}
If $u$ is an abstract function with values in ${\dot{H}^q}$ ($q\in \mathbb{R}$), then the integrals which appears in the above definitions are taken in Bochner's sense. A measurable function $u : [0, \infty)\rightarrow {\dot{H}^q}$ is Bochner integrable if $|u|_q$ is Lebesgue integrable.

Let $(\Omega, \mathcal{F}, \{\mathcal{F}_t \}_{t\geq 0}, {P})$ be a complete filtered probability space satisfying that $\mathcal{F}_0$ contains all $P$-null sets of $\mathcal{F}$. Suppose that $\frac{\partial^2{W}(t,x)}{\partial t \partial x}$ is an infinite dimensional noise defined on $(\Omega, \mathcal{F}, \{\mathcal{F}_t \}_{t\geq 0}, {P})$ such that
\begin{equation}\label{eq1.2}
 \frac{\partial^2{W}(t,x)}{\partial t \partial x} =\sum\limits_{k=1}^{\infty}\sigma_k(t)\dot{\xi}_k(t)e_k(x),
\end{equation}
where $\sigma_k(t)$ is a continuous function, and $\dot{\xi}_k(t)=\frac{d\xi_k(t)}{dt}$, $k=1,2,\ldots$ is the derivative of the standard Wiener process $\xi_k(t),$ $k=1,2,\ldots.$ 

We consider the stochastic space-time fractional wave
equations with $1<\alpha<2$ and $\frac{1}{2}<\beta\leq 1$, as stated in the introduction section,
\begin{equation} \label{eq1.1}
\left\{ \begin{array}
 {l@{\quad} l}
\displaystyle \partial_t^{\alpha}u(t,x)+(-\Delta)^\beta u(t,x)=\frac{\partial^2W(t,x)}{\partial t \partial x},
 \quad 0<t<T,~~ 0<x<1,\\
 \\
  u(t,0)=u(t,1)=0, \quad 0<t<T,\\
  \\
  u(0,x)=v_1(x),~  \partial_tu(0,x)=v_2(x), \quad  0<x<1,\\
 \end{array}
 \right.
\end{equation}
where $T>0$ is a fixed time, $\partial_t^{\alpha}$ denotes the left-sided Caputo fractional derivative of order $\alpha$
with respect to $t$ and $(-\Delta)^\beta$ is the fractional Laplacian, as described above, whose definition is adopted as in \cite{Nochetto}.

Now we define a partition of $[0, T]$ by intervals $[t_i, t_{i+1}]$ for $i = 1, 2, \ldots, I$,
where $t_i = (i-1)\Delta t$, $\Delta t = T/I$. A sequence of noise which approximates the space-time white noise is defined as \cite{du}
\begin{equation}\label{eq1.3}
\displaystyle \frac{\partial^2{W_n}(t,x)}{\partial t \partial x}
 =\sum\limits_{k=1}^{\infty}\sigma^n_k(t)e_k(x)
\left(\sum\limits_{i=1}^{I}\frac{1}{\sqrt{\Delta t}}{\xi}_{ki}\chi_i(t)\right),
\end{equation}
where $\chi_i(t)$ is the characteristic function for the $i$th time subinterval,
\[\xi_{ki} =\frac{1}{\sqrt{\Delta t}}\int_{t_i}^{t_{i+1}}d\xi_k(t)=\frac{1}{\sqrt{\Delta t}}(\xi_k(t_{i+1})-\xi_k(t_i))\sim  N(0, 1),\]
and $\sigma^n_k (t)$ is the approximation of $\sigma_k(t)$ in the space direction. More precisely, replacing $\sigma_k(t)$ by $\sigma^n_k (t)$, we get the noise approximation in space; and replacing $\dot{\xi}_k(t)$ by
$\sum_{i=1}^{I}\frac{1}{\sqrt{\Delta t}}{\xi}_{ki}\chi_i(t)$, we get the noise approximation in time. Then $\frac{\partial^2{W_n}(t,x)}{\partial t \partial x}$ is substituted for $\frac{\partial^2{W}(t,x)}{\partial t \partial x}$ in \eqref{eq1.1} to get the following equation:
\begin{equation} \label{eq1.4}
\left\{ \begin{array}
 {l@{\quad} l}
\displaystyle \partial_t^{\alpha}u_n(t,x)+(-\Delta)^\beta u_n(t,x)=\frac{\partial^2W_n(t,x)}{\partial t \partial x},
 \quad 0<t<T,~~0<x<1,\\
 \\
 u_n(t,0)=u_n(t,1)=0, \quad 0<t<T,\\
 \\
  u_n(0,x)=v_1(x),~ \partial_tu_n(0,x)=v_2(x),\quad  0<x <1.\\
 \end{array}
 \right.
\end{equation}
We define $L^2(\Omega; \dot{H}^q)$ as the separable Hilbert space of all strongly-measurable, square-integrable
 random variables $x$, with values in $\dot{H}^q$, such that
\[\|x\|^2_{L^2(\Omega; \dot{H}^q)}=\mathbf{E}|x|_q^2,\]
where $\mathbf{E}$ denotes the expectation.

The following lemma shows a representation of the mild solution to problem \eqref{eq1.1} using the Dirichlet
eigenpairs $\{(\lambda_k, e_k )\}$.
\begin{lemma}\label{le1}
The solution $u$ to problem \eqref{eq1.1} with $1<\alpha<2$ and $\frac{1}{2}<\beta\leq 1$ is given by
\begin{equation}\label{eq2.7}
\begin{split}
\displaystyle u(t,x)&=\int_0^1\mathcal{T}_{\alpha,\beta}(t,x,y)v_1(y)dy
+\int_0^1\mathcal{R}_{\alpha,\beta}(t,x,y)v_2(y)dy\\
\displaystyle &\relphantom{=}{}
+\int_{0}^t\int_{0}^1\mathcal{S}_{\alpha,\beta}(t-s,x,y)d{W}(s,y).
\end{split}
\end{equation}
Here
\begin{equation}\label{eq2.4}
\mathcal{T}_{\alpha,\beta}(t,x,y)
=\sum\limits_{k=1}^{\infty}E_{\alpha,1}(-\lambda_k^\beta t^\alpha)e_k(x)e_k(y)
\end{equation}
is the fundamental solution of
\[\partial_t^\alpha v(t, x)+ (-\Delta)^\beta v(t, x) = 0,\quad v(t,0)=v(t,1)=0,~~v(0,x)=\phi(x),~~\partial_tv(0, x) = 0,\]
so that $v(t, x) =\int^1_0 \mathcal{T}_{\alpha,\beta}(t,x,y)\phi(y)dy$, and
\begin{equation}\label{eq2.5}
\mathcal{R}_{\alpha,\beta}(t,x,y) =\sum\limits_{k=1}^{\infty}t E_{\alpha,2}(-\lambda_k^\beta t^\alpha)e_k(x)e_k(y)
\end{equation}
is the fundamental solution of
\[\partial_t^\alpha v(t, x)+ (-\Delta)^\beta v(t, x) = 0,\quad v(t,0)=v(t,1)=0,~~v(0,x)=0,~~\partial_tv(0, x) = \psi(x),\]
so that $v(t, x) =\int^1_0 \mathcal{R}_{\alpha,\beta}(t,x,y)\psi(y)dy$. For the equation \eqref{eq1.1} but with initial data $v(0,x)=\partial_tv(0,x) \equiv 0$, we shall use the operator defined by
\begin{equation}\label{eq2.6}
\mathcal{S}_{\alpha,\beta}(t,x,y)=\sum\limits_{k=1}^{\infty}t^{\alpha-1}E_{\alpha,\alpha}(-\lambda_k^\beta t^\alpha)e_k(x) e_k(y),
\end{equation}
and $v(t, x) =\int^t_0 \int^1_0 \mathcal{S}_{\alpha,\beta}(t-s,x,y)dW(s,y).$
\end{lemma}
The proof of Lemma \ref{le1} is given in the Appendix.

Similarly, the integral formulation of \eqref{eq1.4} is
\begin{equation}\label{eq2.8}
\begin{split}
\displaystyle u_n(t,x)&=\int^1_0\mathcal{T}_{\alpha,\beta}(t,x,y)v_1(y)dy
+\int_0^1\mathcal{R}_{\alpha,\beta}(t,x,y)v_2(y)dy\\
&\relphantom{=}{}
+\int_{0}^t\int_0^1\mathcal{S}_{\alpha,\beta}(t-s,x,y)d{W_n}(s,y).
\end{split}
\end{equation}
\subsection{Stability and approximation of the noise.}
In this subsection and in the sequel, the symbol $C$ will denote a generic constant,
whose value may change from one line to another.

First, we give the following stability estimates of the homogeneous problem which will play a key role in the error analysis of the finite element method (FEM) approximations.
\begin{lemma}\label{le2.2}
The solution $u(t)$ to the homogeneous problem of \eqref{eq1.1} satisfies for $t>0$
\begin{equation}\label{eq2.11}
\displaystyle | u(t)|_p \leq \left\{\begin{array}{l}Ct^{\frac{\alpha (q-p)}{2\beta}}|v_1|_q+Ct^{1-\frac{\alpha (p-r)}{2\beta}}|v_2|_r,\quad 0 \leq q, r\leq p\leq 2\beta,\\
\\
Ct^{-\alpha}|v_1|_q+Ct^{1-\alpha}|v_2|_r, \quad\quad\qquad\qquad~~~~~~~~~~~ q, r>p,
\end{array}\right.
\end{equation}
and
\begin{equation*}
\displaystyle |\partial_t^\alpha u(t)|_p \leq Ct^{-\frac{\alpha}{2} (2+\frac{p-q}{\beta})}|v_1|_q+Ct^{1-\frac{\alpha}{2} (2+\frac{p-r}{\beta})}|v_2|_r,~~ 0\leq p \leq q, r\leq p+2\beta.
\end{equation*}
\end{lemma}

The proof of Lemma \ref{le2.2} is given in the Appendix.

Under assumptions on $\{\sigma_k(t)\}$ and $\{\sigma^n_k (t)\}$, our first main result shows that the solution $u_n$ of \eqref{eq1.4} indeed approximates $u$, the solution of \eqref{eq1.1}.
\begin{theorem}\label{thm2.1} Assume that $\{\sigma_k(t)\}$ and its derivative are uniformly bounded by
\[|\sigma_k(t)| \leq \mu_k,~ |\sigma^\prime_k(t)|\leq \gamma_k \quad \forall t \in [0, T],\]
and that the coefficients $\{\sigma^n_k (t)\}$ are constructed such that
\[|\sigma_k(t)- \sigma^n_k (t)| \leq \eta^n_k , ~|\sigma^n_k (t)|\leq \mu^n_k ,
~|(\sigma^{n}_k)^\prime(t)| \leq \gamma^n_k \quad \forall t \in [0, T]\]
with positive sequences $\{\eta^n_k\}$ being arbitrarily chosen, $\{\mu^n_k\}$ and
$\{\gamma^n_k\}$ being related to $\{\eta^n_k \mu_k\}$ and $\{\gamma_k\}$. Let $u_n$ and $u$
be the solutions of \eqref{eq1.4} and \eqref{eq1.1}, respectively. Then for some constant $C > 0$
independent of $\Delta t$,  any $0<\varepsilon<\frac{1}{2}$ and $1<\alpha\leq\frac{3}{2}$, we have
\begin{equation}\label{eq2.9}
\begin{split}
\displaystyle \mathbf{E}\|u(t) - u_n(t)\|^2 &\leq
C\sum\limits_{k=1}^{\infty}\lambda_k^{-\frac{2\beta(\alpha-1)}{\alpha}}(\eta_k^n)^2
+C(\Delta t)^2\sum\limits_{k=1}^{\infty}\lambda_k^{-\frac{2\beta(\alpha-1)}{\alpha}}(\gamma_k^n)^2\\
&\relphantom{=}{}
+Ct^{2\varepsilon}(\Delta t)^{2\alpha-1-2\varepsilon}\sum\limits_{k=1}^{\infty}(\mu_k^n)^2, \quad t>0,
\end{split}
\end{equation}
and for $\frac{3}{2}<\alpha<2$,
\begin{equation}\label{eq2.10}
\begin{split}
\displaystyle \mathbf{E}\|u(t) - u_n(t)\|^2
&\leq C\sum\limits_{k=1}^{\infty}\lambda_k^{-\frac{2\beta(\alpha-1)}{\alpha}}(\eta_k^n)^2
+C(\Delta t)^2\sum\limits_{k=1}^{\infty}\lambda_k^{-\frac{2\beta(\alpha-1)}{\alpha}}(\gamma_k^n)^2\\
&\relphantom{=}{}+Ct^{2\alpha-3}(\Delta t)^2
\sum\limits_{k=1}^{\infty}(\mu_k^n)^2, \quad t>0,
\end{split}
\end{equation}
provided that the infinite series are all convergent.
\end{theorem}
The detailed proof of Theorem \ref{thm2.1} is given in the Appendix.
\section{Galerkin Finite Element Approximation}\label{sec:3}
In this section, we provide the Galerkin FEM scheme and derive the corresponding error estimates.
\subsection{Spatially standard Galerkin FEM and its properties}
For the sake of simplicity, we consider the same partition of $[0, 1]$:
$0 = x_1 < x_2 < \cdots < x_{N+2} = 1$ with $x_j = (j-1)h$ and $h = 1/(N+1)$.
Let $V_h$ be the finite element subspace (with the order of the piecewise polynomial bigger or equal to 1) of $\dot{H}^\beta(0,1)$.

On the space $V_h$ we define the orthogonal $L_2$-projection $P_h: {H} \rightarrow V_h$ and the generalized Ritz projection $R_h : \dot{H}^\beta \rightarrow V_h$, respectively, by
$$(P_h\chi, v)=(\chi, v) \quad \forall v \in V_h,$$ and
$$\left((-\Delta)^{\frac{\beta}{2}} R_h\chi, (-\Delta)^{\frac{\beta}{2}}v\right)
=\left((-\Delta)^{\frac{\beta}{2}}\chi, (-\Delta)^{\frac{\beta}{2}}v\right) \quad \forall v \in V_h.$$
Since $\dot{H}^q$ is equivalent to the fractional Sobolev space $H_0^q$ for $q \in [0,1]$ \cite{Thomee},
we have the
following error estimates for $P_h\psi$ and $R_h\psi$.
\begin{lemma}\label{le3.3} The operators $P_h$ and $R_h$ satisfy
\begin{equation}\label{3.1}
\|P_h\psi-\psi\|+h^{\beta}\|(-\Delta)^{\frac{\beta}{2}}(P_h\psi-\psi)\| \leq Ch^q|\psi|_q \quad {\rm for}\quad \psi\in \dot{H}^q,~ q \in [\beta, 2\beta],
\end{equation}
and
\begin{equation}\label{3.2}
\|R_h\psi- \psi\|+h^{\beta}\|(-\Delta)^{\frac{\beta}{2}}(R_h\psi-\psi)\| \leq Ch^q|\psi|_q \quad{\rm for}\quad \psi\in \dot{H}^q,~ q \in [\beta, 2\beta].
\end{equation}
\end{lemma}

Upon introducing the discrete fractional Laplacian $(-\Delta_h)^\beta : V_h \rightarrow V_h$ defined by
\begin{equation}\label{eq3.0}
((-\Delta_h)^\beta \psi, \chi) = ((-\Delta)^\frac{\beta}{2}\psi, (-\Delta)^\frac{\beta}{2} \chi) \quad \forall \psi, \chi \in V_h,
\end{equation}
we can write the spatial FEM approximation of \eqref{eq1.4} as
\begin{equation}\label{eq3.1}
\partial_t^\alpha u^h_n(t)+(-\Delta_h)^\beta u^h_n(t)=P_h\frac{\partial^2 W_n(t,x)}{\partial t \partial x}, \quad 0<t\leq T,~{\rm with}~u^h_n(0)=v_1^h, ~\partial_t{u^h_n}(0)=v_2^h,
\end{equation}
where $v_1^h=P_h v_1$, $v_2^h=P_h v_2$ or $v_2^h=R_h v_2$.

 Now we give a representation of the solution of \eqref{eq3.1} using the eigenvalues and eigenfunctions $\{\lambda_k^{h,\beta}\}_{k=1}^{N}$ and $\{e_k^h\}_{k=1}^{N}$ of the discrete fractional Laplacian $(-\Delta_h)^\beta$. As we know that the operator $(-\Delta_h)^\beta$ is symmetrical, therefore $\{e_k^h\}_{k=1}^{N}$ is orthogonal. Take $\{e_k^h\}_{k=1}^{N}$ as the orthonormal bases in $V_h$ and define the discrete analogues of \eqref{eq2.4}-\eqref{eq2.6} by
\begin{equation}\label{eq3.2}
\mathcal{T}_{\alpha,\beta}^h(t,x,y)= \sum\limits^N_{k=1} E_{\alpha,1}(-\lambda_k^{h,\beta} t^\alpha)e^h_k(x)e^h_k(y),
\end{equation}
\begin{equation}\label{eq3.6}
\mathcal{R}_{\alpha,\beta}^h(t,x,y)= \sum\limits^N_{k=1} tE_{\alpha,2}(-\lambda_k^{h,\beta} t^\alpha)e^h_k(x)e^h_k(y),
\end{equation}
and
\begin{equation}\label{eq3.3}
\mathcal{S}_{\alpha,\beta}^h(t,x,y)= \sum\limits^N_{k=1} t^{\alpha-1} E_{\alpha,\alpha}(-\lambda_k^{h,\beta} t^\alpha) e^h_k(x)e^h_k(y).
\end{equation}
Then the solution $u^h_n$ of the discrete problem \eqref{eq3.1} can be expressed by
\begin{equation}\label{eq3.4}
\begin{split}
\displaystyle u_n^h(t, x) &=\int_0^1\mathcal{T}_{\alpha,\beta}^h(t,x,y)v_1^h(y)dy
 +\int_0^1\mathcal{R}_{\alpha,\beta}^h(t,x,y)v_2^h(y)dy\\
\displaystyle &\relphantom{=}{} + \int^t_0\int_0^1\mathcal{S}_{\alpha,\beta}^h(t - s,x,y)P_hdW_n(s,y).
  \end{split}
\end{equation}

Also, on the finite element space $V_h$, we introduce the discrete norm $|\cdot|_{p,h}$ for any $p\in \mathbb{R}$ defined by
\begin{equation}\label{eq3.7}
|\psi|_{p,h}^2=\sum\limits^N_{k=1}(\lambda_k^{h,\beta})^{\frac{p}{\beta}}(\psi,e^h_k)^2, \quad \psi\in V_h.
\end{equation}
Clearly, the norm $|\cdot|_{p,h}$ is well defined for all real $p$. By the very definition of the
discrete fractional Laplacian $(-\Delta_h)^\beta$ we have $|\psi|_{p,h}=|\psi|_p$ for $p=0, \beta$ and $\forall \psi\in V_h$. So there is no confusion in using $|\psi|_p$ instead of $|\psi|_{p,h}$ for $p=0, \beta$ and $\forall\psi\in V_h$.
Further, we need the following inverse inequality.
\begin{lemma} For any $l >s$, there exists a constant $C$ independent of $h$ such that
\begin{equation}\label{eq3.5}
\displaystyle |\chi|_{l,h}\leq Ch^{s-l}|\chi|_{s,h} \quad \forall \chi\in V_h.
\end{equation}
\end{lemma}
\vskip 5pt\noindent {\bf Proof.} For all $\chi\in V_h$, the inverse inequality $|\chi|_{\beta} \leq Ch^{-\beta}\|\chi\|$
holds \cite[Lemma 4.5.3]{Brenner}. By the definition of $(-\Delta_h)^\beta$, there exists $\max_{1\leq j\leq N} \lambda^{h,\beta}_j \leq Ch^{-2\beta}$.
Thus, for the norm $|\cdot|_{p,h}$ defined in \eqref{eq3.7}, there holds for any real $l > s$,
\[|\chi|_{l,h}^2 \leq C\max_{1\leq j\leq N}(\lambda^{h,\beta}_j)^{\frac{l-s}{\beta}}\sum_{j=1}^N (\lambda^{h,\beta}_j)^{\frac{s}{\beta}} (\chi,e_k^h)^2 \leq C h^{2(s-l)} |\chi|_{s,h}^2.\qquad \Box\]
The following estimates are crucial for the error analysis in what follows.
\begin{lemma}\label{le3.1}
Let $\mathcal{T}_{\alpha,\beta}^h(t,x,y)$ be defined by \eqref{eq3.2} and $v_1^h \in V_h$. Then for all $t>0$,
\begin{equation*}
\displaystyle \bigg|\int_0^1\mathcal{T}_{\alpha,\beta}^h(t,x,y)v_1^h(y)dy\bigg|_{p,h} \leq \left\{\begin{array}{l}Ct^{\frac{\alpha (q-p)}{2\beta}}|v_1^h|_{q,h},
\quad 0\leq q \leq p \leq 2\beta,\\
\\
Ct^{-\alpha}|v_1^h|_{q,h}, \quad\quad\quad~~ q>p.
\end{array}\right.
\end{equation*}
\end{lemma}
\begin{lemma}\label{le3.4}
Let $\mathcal{R}_{\alpha,\beta}^h(t,x,y)$ be defined by \eqref{eq3.6} and $v_2^h \in V_h$. Then for all $t>0$,
\begin{equation*}
\displaystyle \bigg|\int_0^1\mathcal{R}_{\alpha,\beta}^h(t,x,y)v_2^h(y)dy\bigg|_{p,h} \leq \left\{\begin{array}{l}Ct^{1-\frac{\alpha (p-q)}{2\beta}}|v_2^h|_{q,h},
\quad 0 \leq q\leq p\leq 2\beta,\\ \\
Ct^{1-\alpha}|v_2^h|_{q,h}, \quad\quad\quad~~ q>p.
\end{array}\right.
\end{equation*}
\end{lemma}
The proofs of Lemmas \ref{le3.1} and \ref{le3.4} are discrete analogues of those formulated in \eqref{eq2.11}, so here we omit them.
\begin{lemma}\label{le3.2}
Let $\mathcal{S}_{\alpha,\beta}^h(t,x,y)$ be defined by \eqref{eq3.3} and $\psi\in V_h$. Then for all $t > 0$,
\begin{equation*}
\displaystyle \bigg|\int_0^1\mathcal{S}_{\alpha,\beta}^h(t,x,y)\psi(y)dy\bigg|_{p,h} \leq
\left\{\begin{array}{l} Ct^{-1+\alpha+\frac{\alpha(q-p)}{2\beta}}|\psi|_{q,h}, ~~~  p-2\beta \leq q \leq p,\\
\\
\displaystyle Ct^{-1} |\psi|_{q,h},~~~~\quad\quad\quad  q>p.
\end{array}
\right.
\end{equation*}
\end{lemma}
The proof of Lemma \ref{le3.2} is given in the Appendix.
\subsection{Regularity of the solution of \eqref{eq1.4}.}
 We have the following theorem for the regularity of the solution of \eqref{eq1.4}.
\begin{theorem}\label{thm3.1}
Let $u_n$ be the solution of \eqref{eq1.4}. Assume that $\{\sigma^n_k (t)\}$ are uniformly bounded by $|\sigma^n_k(t)|\leq \mu_k^n$ $\forall t\in[0,T]$. Further assume that $v_1 \in L^2(\Omega;\dot{H}^{q})$, $v_2 \in L^2(\Omega; \dot{H}^{r})$, $q, r\in[0,2\beta]$. Then for some constant $C>0$ independent of $\Delta t$, we have for all $t>0$,
\begin{equation}\label{3.6}
\displaystyle \mathbf{E}|u_n(t)|_{2\beta}^2 \leq C t^{\frac{\alpha(q-2\beta)}{\beta}}\mathbf{E}|v_1|_q^2
+Ct^{2-\frac{\alpha(2\beta-r)}{\beta}}\mathbf{E}|v_2|_r^2+Ct^\alpha\frac{1}{\Delta t}
\sum\limits_{k=1}^{\infty}\lambda_k^{\frac{\beta(\alpha+1)}{\alpha}}(\mu^n_k)^2,
\end{equation}
and
\begin{equation}\label{3.7}
\begin{split}
\mathbf{E}\|\partial_t^\alpha u_n(t)\|^2 &\leq Ct^{\frac{\alpha(q-2\beta)}{\beta}}\mathbf{E}|v_1|_q^2+ Ct^{2-\frac{\alpha(2\beta-r)}{\beta}}\mathbf{E}|v_2|_r^2+Ct^\alpha\frac{1}{\Delta t} \sum\limits_{k=1}^{\infty}\lambda_k^{\frac{\beta(\alpha+1)}{\alpha}}(\mu_k^n)^2\\
&\relphantom{======}{}
+C\frac{1}{{\Delta t}}\sum\limits_{k=1}^{\infty}(\mu_k^n)^2,
\end{split}
\end{equation}
provided that the infinite series are all convergent.
\end{theorem}
\vskip 5pt\noindent {\bf Proof.}
It follows from Lemmas \ref{le1} and \ref{le2.2} that
\begin{equation}\label{3.4}
\begin{split}
\displaystyle\mathbf{E}|u_n(t)|_{2\beta}^2
&\leq C\mathbf{E}\bigg|\int_0^1\mathcal{T}_{\alpha,\beta}(t,x,y)v_1(y)dy\bigg|_{2\beta}^2
+C\mathbf{E}\bigg|\int_0^1\mathcal{R}_{\alpha,\beta}(t,x,y)v_2(y)dy\bigg|_{2\beta}^2\\
&  \relphantom{=}{}
+C\mathbf{E}\bigg|\int_0^t\int_0^1\mathcal{S}_{\alpha,\beta}(t-s,x,y)dW_n(s,y)\bigg|_{2\beta}^2\\
& \leq Ct^{\frac{\alpha(q-2\beta)}{\beta}}\mathbf{E}|v_1|_{q}^2+Ct^{2-\frac{\alpha (2\beta-r)}{\beta}}\mathbf{E}|v_2|_{r}^2\\
& \relphantom{=}{}
+C\mathbf{E}\bigg|\int_0^t(t-s)^{\alpha-1}\sum\limits_{k=1}^{\infty}\int_0^1
E_{\alpha,\alpha}(-\lambda_k^\beta(t-s)^\alpha)e_k(y)e_k(x)dW_n(s,y)\bigg|_{2\beta}^2.
\end{split}
\end{equation}
Let
\[\mathcal{F}_1=C\mathbf{E}\bigg|\int_0^t(t-s)^{\alpha-1}\sum\limits_{k=1}^{\infty}\int_0^1
E_{\alpha,\alpha}(-\lambda_k^\beta(t-s)^\alpha)e_k(y)e_k(x)dW_n(s,y)\bigg|_{2\beta}^2.\]
Without loss of generality, we assume that there exists a positive integer $I_t$ such that $t = t_{I_t+1} $. Since $\{e_k\}_{k=1}^{\infty}$ is an orthonormal basis in $L^2(0,1)$, then by H\"{o}lder's inequality and the boundedness assumption on $\sigma^n_k(t)$, we have
\begin{equation}\label{3.5}
\begin{split}
\mathcal{F}_1&=C\mathbf{E}\sum\limits_{k=1}^{\infty}\lambda_k^{2\beta}\left(\int_0^t(t-s)^{\alpha-1}
\sum\limits_{l=1}^{\infty}\int_0^1 E_{\alpha,\alpha}(-\lambda_l^\beta(t-s)^\alpha)e_l(y)e_l(x)dW_n(s,y),e_k\right)^2\\
&=C\mathbf{E}\sum\limits_{k=1}^{\infty}\lambda_k^{2\beta}\left(\int_0^t(t-s)^{\alpha-1}
\int_0^1 E_{\alpha,\alpha}(-\lambda_k^\beta(t-s)^\alpha)e_k(y)dW_n(s,y)\right)^2\\
&=C\sum\limits_{k=1}^{\infty}\lambda_k^{2\beta}\mathbf{E}\Bigg(\int_0^t(t-s)^{\alpha-1}
\int_0^1 E_{\alpha,\alpha}(-\lambda_k^\beta(t-s)^\alpha)e_k(y)\\
\displaystyle & \relphantom{=================}{}\times \sum\limits_{j=1}^{\infty}
\sigma_j^n(s)e_j(y)\left(\sum\limits_{i=1}^{I_t}\frac{1}{\sqrt{\Delta t}}\xi_{ji}\chi_i(s)\right)dyds\Bigg)^2\\
\displaystyle &=C\sum\limits_{k=1}^{\infty}\lambda_k^{2\beta}\mathbf{E}\Bigg(\int_0^t(t-s)^{\alpha-1}
E_{\alpha,\alpha}(-\lambda_k^\beta(t-s)^\alpha)\sigma_k^n(s)\left(\sum\limits_{i=1}^{I_t}\frac{1}{\sqrt{\Delta t}}\xi_{ki}\chi_i(s)\right)ds\Bigg)^2\\
\displaystyle & \leq C\sum\limits_{k=1}^{\infty}\lambda_k^{2\beta}\int_0^t(t-s)^{\alpha-1}ds\mathbf{E}\int_0^t(t-s)^{\alpha-1}\Bigg(
E_{\alpha,\alpha}(-\lambda_k^\beta(t-s)^\alpha)\\
\displaystyle & \relphantom{========================}{}\times
\sigma_k^n(s)\left(\sum\limits_{i=1}^{I_t}\frac{1}{\sqrt{\Delta t}}\xi_{ki}\chi_i(s)\right)\Bigg)^2ds\\
\displaystyle & = C\sum\limits_{k=1}^{\infty}\lambda_k^{2\beta}t^\alpha
\Bigg(\sum\limits_{l=1}^{I_t}\frac{1}{(\Delta t)^2}
\int_{t_l}^{t_{l+1}}(t-s)^{\alpha-1}E^2_{\alpha,\alpha}(-\lambda_k^\beta(t-s)^\alpha)\\
\displaystyle & \relphantom{======================}{}\times
(\sigma_k^n(s))^2\mathbf{E}\left(\xi_{k}(t_{l+1})-\xi_{k}(t_{l})\right)^2ds\Bigg)\\
\displaystyle & \leq Ct^\alpha \frac{1}{\Delta t}\sum\limits_{k=1}^{\infty}\lambda_k^{2\beta}(\mu_k^n)^2
\int_{0}^{t}(t-s)^{\alpha-1}E^2_{\alpha,\alpha}(-\lambda_k^\beta(t-s)^\alpha)ds\\
\displaystyle & \leq Ct^\alpha \frac{1}{\Delta t}\sum\limits_{k=1}^{\infty}\lambda_k^{2\beta}(\mu_k^n)^2
\lambda_k^{-\frac{\beta(\alpha-1)}{\alpha}}\int_{0}^{t}
\bigg|\frac{(\lambda_k^\beta(t-s)^{\alpha})^{\frac{\alpha-1}{2\alpha}}}
{1+\lambda_k^\beta(t-s)^\alpha}\bigg|^2ds\\
\displaystyle & \leq C t^{\alpha}\frac{1}{\Delta t}\sum\limits_{k=1}^{\infty}
\lambda_k^{\frac{\beta(\alpha+1)}{\alpha}}(\mu_k^n)^2,
\end{split}
\end{equation}
where we have used $\frac{(\lambda_k^\beta(t-s)^{\alpha})^{\frac{\alpha-1}{2\alpha}}}{1+\lambda_k^\beta(t-s)^\alpha}\leq C$.

Thus, the conclusion \eqref{3.6} follows immediately from \eqref{3.4} and \eqref{3.5}.

Now we show \eqref{3.7}. Note that $\{e_k\}_{k=1}^{\infty}$ is an orthonormal basis in $L^2(0,1)$, and for any $t\in[0,T]$, there exists $i\in\{1,2,\ldots,I\}$ such that $t\in[t_i, t_{i+1}]$. By \eqref{eq1.3}, \eqref{3.4} and \eqref{3.5}, we deduce that
\begin{equation}\label{3.8}
\begin{split}
\displaystyle\mathbf{E}\|\partial_t^\alpha u_n(t)\|^2
&=\mathbf{E}\int_0^1(\partial_t^\alpha u_n(t,x))^2dx\\
&=\mathbf{E}\int_0^1\left(-(-\Delta)^\beta u_n(t,x)+\sum\limits_{k=1}^{\infty}
\sigma_k^n(t)e_k(x)\left(\sum\limits_{i=1}^{I}\frac{1}{\sqrt{\Delta t}}\xi_{ki}\chi_i(t)\right)\right)^2dx\\
&\leq C\mathbf{E}| u_n(t)|_{2\beta}^2+
C\mathbf{E}\int_0^1\left(\sum\limits_{k=1}^{\infty}
\sigma_k^n(t)e_k(x)\frac{1}{{\Delta t}}(\xi_{k}(t_{i+1})-\xi_{k}(t_{i}))\right)^2dx\\
& \leq C t^{\frac{\alpha(q-2\beta)}{\beta}}\mathbf{E}|v_1|_{q}^2+Ct^{2-\frac{\alpha(2\beta-r)}{\beta}} \mathbf{E}|v_2|_{r}^2+C t^{\alpha}\frac{1}{\Delta t}\sum\limits_{k=1}^{\infty}\lambda_k^{\frac{\beta(\alpha+1)}{\alpha}}
(\mu_k^n)^2\\
&\relphantom{=}{}+C\sum\limits_{k=1}^{\infty}(\sigma_k^n(t))^2\frac{1}{{(\Delta t)^2}}\mathbf{E}(\xi_{k}(t_{i+1})-\xi_{k}(t_{i}))^2\\
& \leq C t^{\frac{\alpha(q-2\beta)}{\beta}}\mathbf{E}|v_1|_{q}^2+C t^{2-\frac{\alpha(2\beta-r)}{\beta}} \mathbf{E}|v_2|_{r}^2+
C t^{\alpha}\frac{1}{\Delta t} \sum\limits_{k=1}^{\infty}\lambda_k^{\frac{\beta(\alpha+1)}{\alpha}}(\mu_k^n)^2\\
&\relphantom{=}{}+C\frac{1}{\Delta t}\sum\limits_{k=1}^{\infty}(\mu_k^n)^2.
\end{split}
\end{equation}

The proof is completed. $\Box$
\subsection{Error estimates for the homogeneous problem.}
First,  we establish error estimates for the homogeneous problem \eqref{eq1.4} with initial data $v_1\in \dot{H}^q, q\in[\beta, 2\beta]$, $v_2=0$.
\begin{theorem}\label{le4.1} Let $u_n$ be the solution of the homogeneous problem \eqref{eq1.4} with $v_1\in \dot{H}^{q}, q\in[\beta, 2\beta]$, $v_2=0$, and $u_n^h$ be the solution of the homogeneous problem \eqref{eq3.1} with $v_1^h=P_hv_1$, $v_2^h=0$. Then with $\ell_h=|\ln h|$,
\begin{equation*}\begin{split}
&\|u_n^h(t)-u_n(t)\|+h^{\beta}\|(-\Delta)^{\frac{\beta}{2}}(u_n^h(t)-u_n(t))\| \leq C \ell_h h^{2\beta}t^{\frac{\alpha q}{2\beta}-\alpha}|v_1|_{q}~~~ \forall t\in(0,T].
\end{split}\end{equation*}
\end{theorem}
\vskip 5pt\noindent {\bf Proof.} For $v_1\in\dot{H}^{q}, q\in[\beta, 2\beta], v_2=0$, we split the error $u_n^h-u_n$ into $$u_n^h-u_n=(u_n^h-P_hu_n)+(P_hu_n-u_n):=\vartheta+\varrho.$$
It follows from Lemma \ref{le2.2} and \eqref{3.1} that
\begin{equation}\label{3.11}\displaystyle \|\varrho(t)\|+h^{\beta}\|(-\Delta)^{\frac{\beta}{2}}\varrho(t)\|
\leq Ch^{2\beta}|u_n(t)|_{2\beta}\leq Ch^{2\beta}t^{\frac{\alpha q}{2\beta}-\alpha}|v_1|_{q}.
\end{equation}
Using the identity $(-\Delta_h)^\beta R_h=P_h(-\Delta)^\beta$, we see that $\vartheta$ satisfies
\[\partial^\alpha_t\vartheta+(-\Delta_h)^\beta\vartheta=(-\Delta_h)^\beta(R_hu_n-P_hu_n)\]
with $\vartheta(0)=\partial_t\vartheta(0)=0$. Then \eqref{eq3.4} implies that
\begin{equation*}\begin{split}
\vartheta(t,x)=
\int^t_0\int_0^1\mathcal{S}_{\alpha,\beta}^h(t-s,x,y)(-\Delta_h)^\beta(R_hu_n(s,y)-P_hu_n(s,y))dyds.
\end{split}\end{equation*}
For any $0<\varepsilon<2\beta$, by \eqref{eq3.5} and Lemmas \ref{le2.2}, \ref{le3.3} and \ref{le3.2},
 we obtain that
\begin{equation}\label{3.12}
\begin{split}
\displaystyle \|\vartheta(t)\|&\leq \int^t_0\bigg{\|} \int_0^1\mathcal{S}_{\alpha,\beta}^h(t-s,x,y)
(-\Delta_h)^\beta(R_hu_n(s,y)-P_hu_n(s,y))dy\bigg{\|}ds\\
\displaystyle & \leq C\int^t_0(t-s)^{\frac{\alpha\varepsilon}{2\beta}-1}
|(-\Delta_h)^\beta(R_hu_n-P_hu_n)(s)|_{\varepsilon-2\beta,h}ds\\
\displaystyle &=C\int^t_0(t-s)^{\frac{\alpha\varepsilon}{2\beta}-1}|(R_hu_n-P_hu_n)(s)|_{\varepsilon,h}ds\\
\displaystyle &\leq Ch^{-\varepsilon}\int^t_0(t-s)^{\frac{\alpha\varepsilon}{2\beta}-1}\|(R_hu_n-P_hu_n)(s)\|ds\\
\displaystyle &\leq
C h^{2\beta-\varepsilon}\int^t_0(t-s)^{\frac{\alpha\varepsilon}{2\beta}-1}
|u_n(s)|_{2\beta}ds\\
\displaystyle &\leq C h^{2\beta-\varepsilon}
\int^t_0(t-s)^{\frac{\alpha\varepsilon}{2\beta}-1}s^{\frac{\alpha q}{2\beta}-\alpha}|v_1|_{q}ds\\
&\leq C B\left(\frac{\alpha\varepsilon}{2\beta}, 1+\frac{\alpha q}{2\beta}-\alpha\right)h^{2\beta-\varepsilon}t^{\frac{\alpha q}{2\beta}-\alpha}|v_1|_{q}\\
&\leq C \varepsilon^{-1}h^{2\beta-\varepsilon}t^{\frac{\alpha q}{2\beta}-\alpha}|v_1|_{q}.
\end{split}
\end{equation}
In the similar way, we have
\begin{equation}\label{3.1611}
\begin{split}
\displaystyle h^\beta\|(-\Delta)^{\frac{\beta}{2}}\vartheta(t)\|&\leq
h^\beta\int^t_0\bigg{|} \int_0^1\mathcal{S}_{\alpha,\beta}^h(t-s,x,y)
(-\Delta_h)^\beta(R_hu_n(s,y)-P_hu_n(s,y))dy\bigg{|}_\beta ds\\
\displaystyle & \leq Ch^\beta\int^t_0(t-s)^{\frac{\alpha\varepsilon}{2\beta}-1}
|(-\Delta_h)^\beta(R_hu_n-P_hu_n)(s)|_{\varepsilon-\beta,h}ds\\
\displaystyle &= Ch^\beta\int^t_0(t-s)^{\frac{\alpha\varepsilon}{2\beta}-1}|(R_hu_n-P_hu_n)(s)|_{\varepsilon+\beta,h}ds\\
\displaystyle &\leq
Ch^{2\beta-\varepsilon}\int^t_0(t-s)^{\frac{\alpha\varepsilon}{2\beta}-1}|(R_hu_n-P_hu_n)(s)|_\beta ds\\
\displaystyle &\leq C h^{2\beta-\varepsilon}\int^t_0(t-s)^{\frac{\alpha\varepsilon}{2\beta}-1}
|u_n(s)|_{2\beta}ds\\
\displaystyle &\leq C h^{2\beta-\varepsilon}
\int^t_0(t-s)^{\frac{\alpha\varepsilon}{2\beta}-1}s^{\frac{\alpha q}{2\beta}-\alpha}|v_1|_{q}ds\\
&\leq C B\left(\frac{\alpha\varepsilon}{2\beta}, 1+\frac{\alpha q}{2\beta}-\alpha\right)h^{2\beta-\varepsilon}t^{\frac{\alpha q}{2\beta}-\alpha}|v_1|_{q}\\
&\leq C \varepsilon^{-1}h^{2\beta-\varepsilon}t^{\frac{\alpha q}{2\beta}-\alpha}|v_1|_{q}.
\end{split}
\end{equation}
The last inequalities in \eqref{3.12} and \eqref{3.1611} follow from the fact $B(\frac{\alpha\varepsilon}{2\beta},1+\frac{\alpha q}{2\beta}-\alpha)
=\frac{\Gamma(\frac{\alpha\varepsilon}{2\beta})\Gamma(1+\frac{\alpha q}{2\beta}-\alpha)}
{\Gamma(\frac{\alpha\varepsilon}{2\beta}+1+\frac{\alpha q}{2\beta}-\alpha)}$ and $\Gamma(\frac{\alpha\varepsilon}{2\beta})\sim \frac{2\beta}{\alpha\varepsilon}$ as $\varepsilon \rightarrow 0^+$, e.g., by means of Laurent expansion of the Gamma function. Then the desired assertion follows by choosing $\varepsilon=1/\ell_h$ and the triangle inequality. $\Box$

Next we  state an  error estimate for the homogeneous problem \eqref{eq1.4} with initial data $v_1=0$, $v_2\in \dot{H}^{2\beta}$.

\begin{theorem}\label{le4.4} Let $u_n$ be the solution of the homogeneous problem \eqref{eq1.4} with $v_1=0, v_2\in \dot{H}^{2\beta}$, and $u_n^h$ be the solution of the homogeneous problem \eqref{eq3.1} with $v_1^h=0$, $v_2^h=R_hv_2$.
Then
\begin{equation*}\begin{split}
&\|u_n^h(t)-u_n(t)\|+h^{\beta}\|(-\Delta)^{\frac{\beta}{2}}(u_n^h(t)-u_n(t))\| \leq C h^{2\beta}t|v_2|_{2\beta}~~~ \forall t\in(0,T].
\end{split}\end{equation*}
\end{theorem}
\vskip 5pt\noindent {\bf Proof.} For $v_1=0, v_2\in \dot{H}^{2\beta}$, we split the error $u_n^h-u_n$ into two terms as $$u_n^h-u_n=(u_n^h-R_hu_n)+(R_hu_n-u_n):=\bar{\vartheta}+\bar{\varrho}.$$
By \eqref{3.2} and Lemma \ref{le2.2}, we have for any $t>0$
\begin{equation}\label{3.13}
\|\bar{\varrho}(t)\|+h^{\beta}\|(-\Delta)^{\frac{\beta}{2}}\bar{\varrho}(t)\| \leq Ch^{2\beta}|u_n(t)|_{2\beta}\leq Ch^{2\beta}t|v_2|_{2\beta}.
\end{equation}
Using the identity $(-\Delta_h)^\beta R_h=P_h(-\Delta)^\beta$, we note that $\bar{\vartheta}$ satisfies
\[\partial_t^\alpha \bar{\vartheta}(t) +(-\Delta_h)^\beta \bar{\vartheta}(t) =- P_h\partial_t^\alpha\bar{\varrho}(t)\quad {\rm with}\quad \bar{\vartheta}(0)=\partial_t\bar{\vartheta}(0)=0.\]
By \eqref{eq3.4}, we obtain that
\[\bar{\vartheta}(t,x)=-\int^t_0 \int^1_0 \mathcal{S}_{\alpha,\beta}^h(t-s,x,y)P_h\partial_s^\alpha\bar{\varrho}(s,y)dyds.\]
Applying Lemmas \ref{le2.2}, \ref{le3.3} and \ref{le3.2}, then we deduce that for $p=0, \beta$,
\begin{equation}\label{3.14}
\begin{split}
\displaystyle |\bar{\vartheta}(t)|_p &\leq \int^t_0\bigg{|}\int^1_0\mathcal{S}_{\alpha,\beta}^h(t-s,x,y)P_h\partial_s^\alpha\bar{\varrho}(s,y)dy\bigg{|}_pds\\
\displaystyle & \leq C \int^t_0(t-s)^{\alpha-1}|\partial_s^\alpha\bar{\varrho}(s)|_pds\\
\displaystyle & =C \int^t_0(t-s)^{\alpha-1}
|(R_h\partial_s^\alpha u_n-\partial_s^\alpha u_n)(s)|_pds\\
\displaystyle & \leq C h^{2\beta-p} \int^t_0(t-s)^{\alpha-1}|\partial_s^\alpha u_n(s)|_{2\beta}ds\\
\displaystyle & \leq C h^{2\beta-p} \int^t_0(t-s)^{\alpha-1}s^{1-\alpha}|v_2|_{2\beta}ds\\
\displaystyle & \leq C h^{2\beta-p} B\left(\alpha,2-\alpha\right)t|v_2|_{2\beta}.
\end{split}
\end{equation}

Thus the conclusion follows immediately by the triangle inequality and \eqref{3.13}-\eqref{3.14}. $\Box$

Finally we show an error estimate for the homogeneous problem \eqref{eq1.4} with initial data $v_1=0$, $v_2\in H$.
\begin{theorem}\label{le4.3} Let $u_n$ be the solution of the homogeneous problem \eqref{eq1.4} with $v_1=0, v_2\in H$, and $u_n^h$ be the solution of the homogeneous problem \eqref{eq3.1} with $v_1^h=0$, $v_2^h=P_hv_2$. Then with $\ell_h=|\ln h|$,
\begin{equation}\label{3.16}
\begin{split}
\|u_n^h(t)-u_n(t)\|+h^{\beta}\|(-\Delta)^{\frac{\beta}{2}}(u_n^h(t)-u_n(t))\|
\leq C\ell_h h^{2\beta}t^{1-\alpha}\|v_2\|\quad \forall t\in(0,T].
\end{split}\end{equation}
\end{theorem}
\vskip 5pt\noindent {\bf Proof.} For $v_1=0, v_2\in H$, we split the error $u_n^h-u_n$ into two terms: $$u_n^h-u_n=(u_n^h-P_hu_n)+(P_hu_n-u_n):=\hat{\vartheta}+\hat{\varrho}.$$
It follows from Lemma \ref{le2.2} and \eqref{3.1} that
\begin{equation*}\displaystyle \|\hat{\varrho}(t)\|+h^{\beta}\|(-\Delta)^{\frac{\beta}{2}}\hat{\varrho}(t)\|
\leq Ch^{2\beta}|u_n(t)|_{2\beta}\leq Ch^{2\beta}t^{1-\alpha}\|v_2\|.
\end{equation*}
Using the identity $(-\Delta_h)^\beta R_h=P_h(-\Delta)^\beta$, we see that $\hat{\vartheta}$ satisfies
\[\partial^\alpha_t\hat{\vartheta}+(-\Delta_h)^\beta\hat{\vartheta}=(-\Delta_h)^\beta(R_hu_n-P_hu_n)\]
with $\hat{\vartheta}(0)=\partial_t\hat{\vartheta}(0)=0$. Then \eqref{eq3.4} implies that
\[\hat{\vartheta}(t,x)=\int^t_0\int_0^1\mathcal{S}_{\alpha,\beta}^h(t-s,x,y)(-\Delta_h)^\beta(R_hu_n(s,y)-P_hu_n(s,y))dyds.\]
For any $0<\varepsilon<2\beta$, by \eqref{eq3.5} and Lemmas \ref{le2.2}, \ref{le3.3} and \ref{le3.2}, we obtain that
for $p=0, \beta$,
\begin{equation}\label{3.2111}
\begin{split}
\displaystyle |\hat{\vartheta}(t)|_p&\leq \int^t_0\bigg{|} \int_0^1\mathcal{S}_{\alpha,\beta}^h(t-s,x,y)
(-\Delta_h)^\beta(R_hu_n(s,y)-P_hu_n(s,y))dy\bigg{|}_pds\\
\displaystyle &\leq C\int^t_0(t-s)^{\frac{\alpha\varepsilon}{2\beta}-1}
|(-\Delta_h)^\beta(R_hu_n-P_hu_n)(s)|_{\varepsilon-2\beta+p,h}ds\\
\displaystyle &= C\int^t_0(t-s)^{\frac{\alpha\varepsilon}{2\beta}-1}
|(R_hu_n-P_hu_n)(s)|_{\varepsilon+p,h}ds\\
\displaystyle &\leq Ch^{-\varepsilon}\int^t_0(t-s)^{\frac{\alpha\varepsilon}{2\beta}-1}
|(R_hu_n-P_hu_n)(s)|_{p}ds\\
\displaystyle &\leq C h^{2\beta-p-\varepsilon}\int^t_0(t-s)^{\frac{\alpha\varepsilon}{2\beta}-1}|u_n(s)|_{2\beta}ds\\
\displaystyle &\leq C h^{2\beta-p-\varepsilon}\int^t_0(t-s)^{\frac{\alpha\varepsilon}{2\beta}-1}
s^{1-\alpha}\|v_2\|ds\\
\displaystyle &\leq CB\left(\frac{\alpha\varepsilon}{2\beta},2-\alpha\right)
h^{2\beta-p-\varepsilon}t^{1-\alpha}\|v_2\|\\
\displaystyle &\leq C\varepsilon^{-1}h^{2\beta-p-\varepsilon}
t^{1-\alpha}\|v_2\|.
\end{split}
\end{equation}
The last inequality follows from the fact $B(\frac{\alpha\varepsilon}{2\beta},2-\alpha)
=\frac{\Gamma(\frac{\alpha\varepsilon}{2\beta})\Gamma(2-\alpha)}
{\Gamma(\frac{\alpha\varepsilon}{2\beta}+2-\alpha)}$ and $\Gamma(\frac{\alpha\varepsilon}{2\beta})\sim \frac{2\beta}{\alpha\varepsilon}$ as $\varepsilon \rightarrow 0^+$. Then the desired assertion follows by choosing $\varepsilon=1/\ell_h$ and the triangle inequality. $\Box$

\begin{remark}\label{re4.4}
Let $u_n$ be the solution of the homogeneous problem \eqref{eq1.4} with $v_1=0, v_2\in \dot{H}^r, r\in[0, 2\beta]$, and $u_n^h$ be the solution of the homogeneous problem \eqref{eq3.1} with $v_1^h=0$, $v_2^h=P_hv_2$. Then with $\ell_h=|\ln h|$,
\begin{equation}\label{3.22}
\begin{split}
\|u_n^h(t)-u_n(t)\|+h^{\beta}\|(-\Delta)^{\frac{\beta}{2}}(u_n^h(t)-u_n(t))\|
\leq C\ell_h h^{2\beta}t^{1-\alpha+\frac{\alpha r}{2\beta}}|v_2|_r\quad \forall t\in(0,T].
\end{split}\end{equation}
\end{remark}

\subsection{Error estimates for the nonhomogeneous problem.}
First, we derive an error estimate for the nonhomogeneous problem \eqref{eq1.4} with initial data $v_1=v_2=0$.
\begin{theorem}\label{le4.2}
Let $u_n$ and $u_n^h$ be the solutions of \eqref{eq1.4} and \eqref{eq3.1} with $v_1=v_2=0$, respectively. Assume that $\{\sigma^n_k (t)\}$ are uniformly bounded by $|\sigma^n_k(t)|\leq \mu_k^n$ ~$\forall t\in[0,T]$.
Then with $\ell_h=|\ln h|$,
\begin{equation*}
\begin{split}
\displaystyle &\mathbf{E}\|u_n^h(t)-u_n(t)\|^2+h^{2\beta}
\mathbf{E}\|(-\Delta)^{\frac{\beta}{2}}(u_n^h(t)-u_n(t))\|^2\\
&\relphantom{=}{}\leq C\ell_h h^{4\beta}\frac{1}{\Delta t}
\sum\limits_{k=1}^{\infty}\lambda_k^{\frac{\beta(\alpha+1)}{\alpha}}(\mu_k^n)^2\quad \forall t\in[0,T],
\end{split}
\end{equation*}
provided that the infinite series are  convergent, where $C$ is a positive constant independent of $\Delta t$ and $h$.
\end{theorem}
\vskip 5pt\noindent {\bf Proof.} We split the error $u_n^h-u_n$ into
$$u_n^h-u_n=(u_n^h-P_hu_n)+(P_hu_n-u_n):=\upsilon+\rho.$$
By \eqref{3.1} we have
\begin{equation}\label{3.9}
\begin{split}
\displaystyle \mathbf{E}\|\rho(t)\|^2+h^{2\beta}\mathbf{E}\|(-\Delta)^{\frac{\beta}{2}}\rho(t)\|^2  &\leq Ch^{4\beta}\mathbf{E}|u_n(t)|_{2\beta}^2.
\end{split}
\end{equation}
Note that $v_1=v_2=0$, hence it follows from \eqref{3.4} and \eqref{3.5} that
\begin{equation}\label{3.15}
\begin{split}
\displaystyle \displaystyle \mathbf{E}|u_n(t)|_{2\beta}^2
&= \mathbf{E}\bigg|\int_0^t\int_0^1\mathcal{S}_{\alpha,\beta}(t-s,x,y)dW_n(s,y)\bigg|_{2\beta}^2\\
&\leq  Ct^{\alpha}\frac{1}{\Delta t}\sum\limits_{k=1}^{\infty}\lambda_k^{\frac{\beta(\alpha+1)}{\alpha}}(\mu_k^n)^2.
\end{split}
\end{equation}
Therefore
\begin{equation}\label{3.3}
\begin{split}
\displaystyle \displaystyle \mathbf{E}\|\rho(t)\|^2+h^{2\beta}\mathbf{E}\|(-\Delta)^{\frac{\beta}{2}}\rho(t)\|^2 &\leq  Ct^{\alpha}h^{4\beta}\frac{1}{\Delta t}\sum\limits_{k=1}^{\infty}\lambda_k^{\frac{\beta(\alpha+1)}{\alpha}}(\mu_k^n)^2.
\end{split}
\end{equation}
Moreover, we consider the equation
\[\partial_t^\alpha \upsilon +(-\Delta_h)^\beta \upsilon = (-\Delta_h)^\beta(R_hu_n-P_hu_n)
\quad{\rm with}\quad\upsilon(0)=\partial_t\upsilon(0)=0.\]
Then it follows from \eqref{eq3.4} that
\[\upsilon(t,x)=\int^t_0 \int_0^1\mathcal{S}_{\alpha,\beta}^h(t-s,x,y)(-\Delta_h)^\beta(R_hu_n(s,y)-P_hu_n(s,y))dyds.\]
For any $0<\varepsilon<2\beta$, by the similar argument of \eqref{3.2111} and using H\"{o}lder's inequality, we deduce that for $p=0, \beta$,
\begin{equation*}
\begin{split}
\displaystyle \mathbf{E}|\upsilon(t)|_p^2 &\leq \mathbf{E}\left(\int^t_0\bigg{|} \int_0^1\mathcal{S}_{\alpha,\beta}^h(t-s,x,y)(-\Delta_h)^\beta(R_hu_n(s,y)-P_hu_n(s,y))dy\bigg{|}_pds\right)^2\\
\displaystyle & \leq C \mathbf{E}\left(\int^t_0(t-s)^{\frac{\alpha\varepsilon}{2\beta}-1}
|(-\Delta_h)^\beta(R_hu_n-P_hu_n)(s)|_{\varepsilon-2\beta+p,h}ds\right)^2\\
\displaystyle & =C \mathbf{E}\left(\int^t_0(t-s)^{\frac{\alpha\varepsilon}{2\beta}-1}
|(R_hu_n-P_hu_n)(s)|_{\varepsilon+p,h}ds\right)^2\\
\displaystyle & \leq C h^{-2\varepsilon}\mathbf{E}\left(\int^t_0(t-s)^{\frac{\alpha\varepsilon}{2\beta}-1}
|(R_hu_n-P_hu_n)(s)|_pds\right)^2\\
\displaystyle & \leq C h^{4\beta-2p-2\varepsilon} \mathbf{E}\left(\int^t_0(t-s)^{\frac{\alpha\varepsilon}{2\beta}-1}
|u_n(s)|_{2\beta}ds\right)^2\\
\displaystyle & \leq C h^{4\beta-2p-2\varepsilon} \int^t_0(t-s)^{\frac{\alpha\varepsilon}{2\beta}-1}ds
\int^t_0(t-s)^{\frac{\alpha\varepsilon}{2\beta}-1}\mathbf{E}|u_n(s)|_{2\beta}^2ds\\
\displaystyle & \leq C h^{4\beta-2p-2\varepsilon} t^{\frac{\alpha\varepsilon}{2\beta}}
\int^t_0(t-s)^{\frac{\alpha\varepsilon}{2\beta}-1}\mathbf{E}|u_n(s)|_{2\beta}^2ds.
\end{split}
\end{equation*}
Using \eqref{3.15}, we obtain for $p=0, \beta$,
\begin{equation}\label{3.10}
\begin{split}
\displaystyle \mathbf{E}|\upsilon(t)|_p^2 &\leq C t^{\frac{\alpha\varepsilon}{2\beta}} h^{4\beta-2p-2\varepsilon} \int^t_0(t-s)^{\frac{\alpha\varepsilon}{2\beta}-1}s^{\alpha} \frac{1}{\Delta t}\sum\limits_{k=1}^{\infty}\lambda_k^{\frac{\beta(\alpha+1)}{\alpha}}(\mu_k^n)^2ds\\
&\leq C t^{\frac{\alpha\varepsilon}{2\beta}}h^{4\beta-2p-2\varepsilon} \frac{1}{\Delta t}
\sum\limits_{k=1}^{\infty}\lambda_k^{\frac{\beta(\alpha+1)}{\alpha}}(\mu_k^n)^2
B\left(\frac{\alpha\varepsilon}{2\beta},1+\alpha\right)\\
&\leq C \varepsilon^{-1}h^{4\beta-2p-2\varepsilon} \frac{1}{\Delta t}
\sum\limits_{k=1}^{\infty}\lambda_k^{\frac{\beta(\alpha+1)}{\alpha}}(\mu_k^n)^2.
\end{split}
\end{equation}
The last inequality follows from the fact $B(\frac{\alpha\varepsilon}{2\beta},1+\alpha)=\frac{\Gamma(\frac{\alpha\varepsilon}{2\beta})\Gamma(1+\alpha)}
{\Gamma(\frac{\alpha\varepsilon}{2\beta}+1+\alpha)}$ and $\Gamma(\frac{\alpha\varepsilon}{2\beta})\sim \frac{2\beta}{\alpha\varepsilon}$ as $\varepsilon\rightarrow 0^+$. Then the desired assertion follows by choosing $\varepsilon=1/\ell_h$ and the triangle inequality. $\Box$

As a simple consequence of Theorem \ref{le4.1}, Remark \ref{re4.4} and Theorem \ref{le4.2}, we have

\begin{theorem}\label{4.1}
Let $u_n$ be the solution of problem \eqref{eq1.4} with $v_1\in\dot{H}^q$, $q\in[\beta,2\beta]$,  $v_2\in\dot{H}^r$, $r\in[0,2\beta]$, and $u_n^h$ be the solution of \eqref{eq3.1} with $v_1^h=P_hv_1, v_2^h=P_hv_2$. Assume that $\{\sigma^n_k (t)\}$ are uniformly bounded by $|\sigma^n_k(t)|\leq \mu_k^n$ $\forall t\in[0,T]$.
Then with $\ell_h=|\ln h|$,
\begin{equation*}\begin{split}
\displaystyle&\mathbf{E}\|u_n(t)-u_n^h(t)\|^2+h^{2\beta}\mathbf{E}\|(-\Delta)^{\frac{\beta}{2}}(u_n^h(t)-u_n(t))\|^2 \leq C \ell_h^2 h^{4\beta}t^{\frac{\alpha q}{\beta}-2\alpha}\mathbf{E}|v_1|_{q}^2\\
&\relphantom{=}{}
+C\ell_h^2 h^{4\beta}t^{2-2\alpha+\frac{\alpha r}{\beta}}\mathbf{E}|v_2|^2_{r}
+C\ell_hh^{4\beta}\frac{1}{\Delta t}\sum\limits_{k=1}^{\infty}
\lambda_k^{\frac{\beta(\alpha+1)}{\alpha}}(\mu^n_k)^2\quad \forall t \in [0, T],
\end{split}
\end{equation*}
provided that the infinite series are  convergent, where $C$ is a positive constant independent of $\Delta t$ and $h$.
\end{theorem}
\vskip 5pt\noindent {\bf Proof.} By \eqref{eq2.8} and \eqref{eq3.4}, we have
\begin{equation*}\begin{split}
\displaystyle u_n(t, x)-u_n^h(t, x) &=\int_0^1\left(\mathcal{T}_{\alpha,\beta}(t,x,y)v_1(y)-\mathcal{T}_{\alpha,\beta}^h(t,x,y)P_hv_1(y)\right)dy\\
\displaystyle&\relphantom{=}{}+\int_0^1\left(\mathcal{R}_{\alpha,\beta}(t,x,y)v_2(y)
-\mathcal{R}_{\alpha,\beta}^h(t,x,y)P_hv_2(y)\right)dy\\
\displaystyle&\relphantom{=}{}+
\int_0^t\int_0^1\left(\mathcal{S}_{\alpha,\beta}(t-s,x,y)-\mathcal{S}_{\alpha,\beta}^h(t-s,x,y)P_h
\right)dW_n(s,y).
\end{split}
\end{equation*}
Then by the triangle inequality, Theorem \ref{le4.1}, Remark \ref{re4.4} and Theorem \ref{le4.2},
we obtain the desired result. $\Box$

Furthermore, thanks to  Theorems \ref{thm2.1} and  \ref{4.1}, a space-time error estimate
for problem  \eqref{eq1.1}  follows from the triangle inequality.
\begin{theorem}\label{4.2}
Suppose that the assumptions of Theorem \ref{thm2.1} hold. Let $u$ be the solution of \eqref{eq1.1} with $v_1\in\dot{H}^q$, $q\in[\beta,2\beta]$, $v_2\in\dot{H}^r$, $r\in[0,2\beta]$, and $u_n^h$ be the solution of \eqref{eq3.1} with $v_1^h=P_hv_1$, $v_2^h=P_hv_2$. Then for any $0<\varepsilon<\frac{1}{2}$ and with $\ell_h=|\ln h|$,
\begin{equation*}
\begin{split}
&\mathbf{E}\|u(t)-u_n^h(t)\|^2 \leq
 C\sum\limits_{k=1}^{\infty}\lambda_k^{-\frac{2\beta(\alpha-1)}{\alpha}}(\eta_k^n)^2
 +C(\Delta t)^2\sum\limits_{k=1}^{\infty}\lambda_k^{-\frac{2\beta(\alpha-1)}{\alpha}}(\gamma_k^n)^2\\
 \displaystyle &\relphantom{=}{}
+C(\Delta t)^{2\alpha-1-2\varepsilon}\sum\limits_{k=1}^{\infty}(\mu_k^n)^2
+C \ell_h^2 h^{4\beta}t^{\frac{\alpha q}{\beta}-2\alpha}\mathbf{E}|v_1|_{q}^2\\
&\relphantom{=}{}
+C\ell_h^2 h^{4\beta}t^{2-2\alpha+\frac{\alpha r}{\beta}}\mathbf{E}|v_2|^2_{r}
+C\ell_hh^{4\beta}\frac{1}{\Delta t}\sum\limits_{k=1}^{\infty}
\lambda_k^{\frac{\beta(\alpha+1)}{\alpha}}(\mu^n_k)^2
\quad{\rm for}\quad 1<\alpha\leq\frac{3}{2},
\end{split}
\end{equation*}
and
\begin{equation*}
\begin{split}
&\mathbf{E}\|u(t)-u_n^h(t)\|^2\leq
 C\sum\limits_{k=1}^{\infty}\lambda_k^{-\frac{2\beta(\alpha-1)}{\alpha}}(\eta_k^n)^2
 +C(\Delta t)^2\sum\limits_{k=1}^{\infty}\lambda_k^{-\frac{2\beta(\alpha-1)}{\alpha}}(\gamma_k^n)^2\\
 \displaystyle &\relphantom{=}{}
+C(\Delta t)^2\sum\limits_{k=1}^{\infty}(\mu_k^n)^2
+C \ell_h^2 h^{4\beta}t^{\frac{\alpha q}{\beta}-2\alpha}\mathbf{E}|v_1|_{q}^2\\
&\relphantom{=}{}
+C\ell_h^2 h^{4\beta}t^{2-2\alpha+\frac{\alpha r}{\beta}}\mathbf{E}|v_2|^2_{r}
+C\ell_hh^{4\beta}\frac{1}{\Delta t}\sum\limits_{k=1}^{\infty}
\lambda_k^{\frac{\beta(\alpha+1)}{\alpha}}(\mu^n_k)^2
\quad{\rm for}\quad \frac{3}{2}< \alpha<2,
\end{split}
\end{equation*}
provided that the infinite series are  all convergent, where $C$ is a positive constant independent of $\Delta t$ and $h$.
\end{theorem}
\section{Numerical Results}\label{sec:4}
Here we present some numerical tests to verify the theoretical error estimates for the Galerkin FEM.
For definiteness, we simulate \eqref{eq1.1}, \eqref{eq1.4} and \eqref{eq3.1} with
\[\sigma_k(t)=\frac{1}{k^3}, \quad
\sigma_k^n(t)=\left\{\begin{array}{ll}\sigma_k(t), & \hbox{$k\leq n$,} \\
\\
0, & \hbox{$k>n$,}\end{array}\right.
\quad v_1(x)=-4x^2+4x {\rm ~~and~~} v_2(x)=x.\]
The upper bounds $\eta^n_k$, $\mu^n_k$ and $\gamma^n_k$ given in Theorem \ref{thm2.1} can be chosen as
\[\eta_k^n(t)=\left\{\begin{array}{ll}0, & \hbox{$k\leq n$,} \\
\\
\frac{1}{k^3}, & \hbox{$k>n$,}\end{array}\right. \qquad \mu^n_k=\gamma^n_k=\frac{1}{k^3}.\]
For getting the numerical errors and convergence orders, according to the definition of Ito's integral, we introduce the reference (``exact") solution $u_{ref}$
defined by
\begin{equation}\label{4.5}
\begin{split}
\displaystyle u_{ref}&=\int_0^1\mathcal{T}_{\alpha,\beta}(t,x,y)v_1(y)dy
+\int_0^1\mathcal{R}_{\alpha,\beta}(t,x,y)v_2(y)dy\\
\displaystyle &\relphantom{=}{} +\sum\limits_{n=1}^{N'}\int_{0}^1\mathcal{S}_{\alpha,\beta}(t-t_n,x,y)\sum\limits_{k=1}^{M'}\sigma_k(t_n)
(\xi_k(t_{n+1})-\xi_k(t_n))e_k(y)dy,
\end{split}
\end{equation}
where the values of the Mittag-Leffler function at each nodal point $t_n$ in \eqref{4.5} are obtained by applying the algorithm given in \cite{Garrappa}; $N^{\prime}=1000$ with the uniform discretization of the time interval; $M^{\prime}=1000$, $\lambda_k=k^2\pi^2$, and $e_k(y)=\sqrt{2}\sin k\pi y$.

In our numerical experiments, we first begin with a study of $L^2$-norm error
in the mean-squared sense between the solutions $u_{ref}$ and  $u_n$, which is based on the
space and time discretization of the noise. By \eqref{A1} in Appendix, we only need to
simulate the stochastic integrals in \eqref{eq2.7} and \eqref{eq2.8}. 
By Theorem \ref{thm2.1}, we need to verify 
\begin{equation}\label{4.6}
(\mathbf{E}\|u_{ref}(t_{I+1})-u_n(t_{I+1})\|^2)^{\frac{1}{2}}\leqslant\left\{\begin{array}{ll}C (\Delta t)^{\alpha-\frac{1}{2}-\varepsilon},
& \hbox{$1<\alpha \leq \frac{3}{2}$,} \\
\\
C \Delta t, & \hbox{$\frac{3}{2}<\alpha<2$.}\end{array}\right.
\end{equation}
Since
\[(\mathbf{E}\|u_{ref}(t_{I+1})-u_n(t_{I+1})\|^2)^{\frac{1}{2}}\approx \left(\frac{1}{M}\sum\limits_{l=1}^{M}\|u_{ref}(t_{I+1},\omega_l)-u_n(t_{I+1},\omega_l)\|^2\right)^{\frac{1}{2}},\]
we take $M=1000$ as the number of simulation trajectories. 
The numerical results with $\beta=0.75$ and $t_{I+1}=T=1$ are presented in Table \ref{table1}. From Table \ref{table1}, it can be observed that, the mean-squared $L^2$-norm errors converge fast as $\Delta t$ decreases, and its convergence rates
are closely related to the values of $\alpha$, which agrees well with \eqref{4.6}.
These results confirm the error estimate in Theorem \ref{thm2.1}.
\begin{table}[htbp] \fontsize{5.0pt}{11pt}\selectfont
\centering
\caption{\label{table1} The mean-squared $L^2$-norm modeling errors and convergence rates of \eqref{eq1.4} with $\beta=0.75$, $t_{I+1}=1$, and $M=1000$.} \vspace{5pt}
\begin{tabular}{|c|c|c|c|c|c|c|c|c|c|c|}
\hline $1/\Delta t$ & $\alpha=1.1$ & Rate & $\alpha=1.25$ & Rate & $\alpha=1.5$ & Rate & $\alpha=1.75$ & Rate & $\alpha=2.0$ & Rate \\
\hline
\multirow{11}{*}{} 25  & 1.2567e-2 & --     & 1.5505e-2 & --     & 1.3168e-2 & --                & 9.8503e-3 & --     & 7.1887e-3 & --    \\
\hline
\multirow{11}{*}{} 50  & 7.6842e-3 & 0.7097 & 9.5036e-3 & 0.7062 & 7.3208e-3 & 0.8470 & 5.1609e-3 & 0.9325 & 3.7898e-3 & 0.9236 \\
\hline
\multirow{11}{*}{} 100 & 4.8080e-3 & 0.6765 & 5.4252e-3 & 0.8088 & 3.7393e-3 & 0.9693 & 2.5460e-3 & 1.0194 & 1.8082e-3 & 1.0676\\
\hline
\multirow{11}{*}{} 125 & 4.2249e-3 & 0.5793 & 4.6975e-3 & 0.6454 & 3.0322e-3 & 0.9392 & 2.0506e-3 & 0.9697 & 1.4707e-3 & 0.9258\\
\hline
\multirow{11}{*}{} 200 & 2.9152e-3 & 0.7895 & 3.1393e-3 & 0.8575 & 2.0008e-3 & 0.8846 & 1.3215e-3 & 0.9348 & 9.4838e-4 & 0.9335 \\
\hline
\end{tabular}
\end{table}

Next, we measure the mean-squared $L^2$-norm error between $u_n$ and the finite element solution $u_n^h$.
We first divide the unit interval $(0, 1)$ into $N + 1$ equally spaced subintervals with
 $h = 1/(N + 1)$. Let $V_h$ be the finite element space consisting of continuous piecewise linear
polynomials. Notice that the eigenpairs $\{\lambda^{h,\beta}_j, e^h_j(x)\}$  of the one-dimensional
discrete fractional Laplacian $(-\triangle_h)^\beta$, defined by \eqref{eq3.0}, satisfies
\begin{equation}\label{4.3}
((-\triangle_h)^\beta e^h_j,\chi)=\lambda^{h,\beta}_j(e^h_j,\chi) \quad\quad \forall \chi\in V_h.
\end{equation}
Then for $j=1, 2, \ldots, N,$
\begin{equation}\label{4.4}
\lambda^{h,\beta}_j= \sum_{k=1}^{\infty}\lambda_k^\beta(e^h_j,e_k)^2.
\end{equation}


The expressions \eqref{4.3}-\eqref{4.4} are used in computing the finite element solution of the Galerkin methods through its representation \eqref{eq3.4}. To validate the convergence rate $2\beta$ for $\|u_n(t_{I+1})-u_n^h(t_{I+1})\|_{L^2(\Omega;H)}$, we simulate the integral expressions in \eqref{eq2.8} and \eqref{eq3.4} with $\Delta t=0.01$ and the number of trajectories $M=500$. The numerical results with $\alpha=1.5, t_{I+1}=T=1$ are presented in Table \ref{table2}. From Table \ref{table2}, we can see that, for fixed $I$, the results of $\|u_n(t_{I+1})-u_n^h(t_{I+1})\|_{L^2(\Omega;H)}$ converge as $h$ approaches to 0, and the corresponding convergence rates at least have an order of $\mathcal{O}(h^{2\beta})$, which in turn justify the statement of Theorem \ref{4.1}.
Combining the datum of the tables \ref{table1} and \ref{table2} gives the convergence rates, at least $\mathcal{O}((\Delta t)^{\alpha-1/2-\varepsilon}+ h^{2\beta})$ for $1<\alpha \leq 3/2$ and
$\mathcal{O}(\Delta t+ h^{2\beta})$ for $3/2<\alpha<2$, where the errors, measured by the mean-squared $L^2$-norm, are the differences between the reference (``exact") solution $u_{ref}(t_{I+1}, \omega_l)$ and the numerical solution $u_n^h(t_{I+1},\omega_l)$. This means that Theorem
\ref{4.2} is valid.
\begin{table}[htbp]
\centering
\caption{\label{table2} The mean-squared $L^2$-norm errors and convergence rates of the FEM approximations \eqref{eq3.1} with $\alpha=1.5$, $t_{I+1}=1$, $\Delta t=0.01$ and $M=500$.} \vspace{5pt}
\begin{tabular}{|c|c|c|c|c|c|c|}
\hline   $1/h$ &$\beta=0.6$  & Rate &$\beta=0.8$  & Rate  &$\beta=1.0$  & Rate\\
\hline
\multirow{7}{*}{} 10  & 7.6377e-3 & --     & 6.8670e-3 & --     & 4.0807e-3 & --     \\
\hline
\multirow{7}{*}{} 25  & 1.3772e-3 & 1.8695 & 1.1233e-3 & 1.9759 & 5.9384e-4 & 2.1035 \\
\hline
\multirow{7}{*}{} 50  & 3.8494e-4 & 1.8390 & 2.7870e-4 & 2.0109 & 1.6592e-4 & 1.8396 \\
\hline
\multirow{7}{*}{} 75  & 1.8444e-4 & 1.8146 & 1.2644e-4 & 1.9493 & 7.3794e-5 & 1.9983 \\
\hline
\multirow{7}{*}{} 100 & 1.0920e-4 & 1.8221 & 6.9840e-5 & 2.0632 & 4.1522e-5 & 1.9989 \\
\hline
\end{tabular}
\end{table}

\section{Conclusions}
Not only the pillars of anomalous diffusion but also the relation between attenuation and frequency of wave propagation in heterogeneous media have power-law structure. Fractional calculus plays a key in characterizing this structure. In general, the stochastic perturbations can not be avoided in physical system; sometimes they even can not be ignored, which need to add the corresponding stochastic term to the deterministic governing equation. The paper considers the stochastic space-time fractional wave equation forced by an additive space-time white noise, which describe the power-law attenuated wave propagation affected by the white Gaussian noise. The effective finite element approximation for the fractional Laplacian is provided, and the definition of the discrete fractional Laplacian is introduced. The special norms and function spaces are established to analyze the regularity and stability of the equations, in particular, the mean-squared $L^2$-norm priori estimates for the finite element approximation error and modeling error, resulted by discretizing the space-time white noise, are derived. The performed numerical simulation results confirm the theoretical analysis.
\section*{Acknowledgments}
This work was supported by the National Natural Science Foundation of China under Grant No. 11271173
and No. 11571153.
\section*{Appendix}\label{sec:5}
\vskip 5pt\noindent {\bf Proof of Lemma \ref{le1}.}
Let
\begin{equation}\label{a1}
\displaystyle u(t,x)= \sum\limits_{k=1}^{\infty}(u(t),e_k)e_k(x)=\sum\limits_{k=1}^{\infty}u_k(t)e_k(x).
\end{equation}
Substituting \eqref{a1} into \eqref{eq1.1}, we get that for $1<\alpha<2$,
\begin{equation}\label{a2}
\displaystyle \partial^\alpha_t u_k(t)+ \lambda_k^\beta u_k(t)= \sigma_k(t)\dot{\xi}_k(t),
\quad u_k(0)=v_{1k},~\partial_tu_k(0)=v_{2k}.
\end{equation}
By Theorem 5.15 in \cite{Kilbas}, we have
\begin{equation}\label{a13}\begin{split}
\displaystyle u_k(t)&= v_{1k} E_{\alpha,1}(-\lambda_k^\beta t^\alpha)+v_{2k}tE_{\alpha,2}(-\lambda_k^\beta t^\alpha)\\
&\relphantom{=}{}+\int^t_0 (t-s)^{\alpha-1}E_{\alpha,\alpha}(-\lambda_k^\beta (t-s)^\alpha)\sigma_k(s)\dot{\xi}_k(s)ds.
\end{split}\end{equation}
Then, it follows from \eqref{a1} and \eqref{a13} that
\begin{equation*}\label{a4}\begin{split}
\displaystyle u(t,x)& = \sum\limits_{k=1}^{\infty}v_{1k} E_{\alpha,1}(-\lambda_k^\beta t^\alpha)e_k(x)
+\sum\limits_{k=1}^{\infty}v_{2k}t E_{\alpha,2}(-\lambda_k^\beta t^\alpha)e_k(x)\\
&\relphantom{=}{}+\sum\limits_{k=1}^{\infty}
\int^t_0 (t-s)^{\alpha-1}E_{\alpha,\alpha}(-\lambda_k^\beta (t-s)^\alpha)\sigma_k(s)\dot{\xi}_k(s)dse_k(x)\\
&=\sum\limits_{k=1}^{\infty} E_{\alpha,1}(-\lambda_k^\beta t^\alpha)(v_1, e_k)e_k(x)
+\sum\limits_{k=1}^{\infty}tE_{\alpha,2}(-\lambda_k^\beta t^\alpha)(v_2, e_k)e_k(x)\\
&\relphantom{=}{}+\sum\limits_{k=1}^{\infty}
\int^t_0 (t-s)^{\alpha-1} E_{\alpha,\alpha}(-\lambda_k^\beta (t-s)^\alpha)
\left(\sum\limits_{l=1}^{\infty}\sigma_l(s)\dot{\xi}_l(s)e_l, e_k\right)e_k(x) ds\\
&=\int_0^1\sum\limits_{k=1}^{\infty} E_{\alpha,1}(-\lambda_k^\beta t^\alpha)e_k(x)e_k(y)v_1(y)dy\\
&\relphantom{=}{}+\int_0^1\sum\limits_{k=1}^{\infty}tE_{\alpha,2}(-\lambda_k^\beta t^\alpha)e_k(x) e_k(y)v_2(y)dy\\
&\relphantom{=}{}
+\int^t_0\int_0^1 (t-s)^{\alpha-1}\sum\limits_{k=1}^{\infty} E_{\alpha,\alpha}(-\lambda_k^\beta (t-s)^\alpha)e_k(x)
e_k(y)d{W}(s,y)\\
&=\int_0^1\mathcal{T}_{\alpha,\beta}(t,x,y)v_1(y)dy
+\int_0^1\mathcal{R}_{\alpha,\beta}(t,x,y)v_2(y)dy\\
&\relphantom{=}{}+\int_{0}^t\int_{0}^1\mathcal{S}_{\alpha,\beta}(t-s,x,y)d{W}(s,y).
\end{split}\end{equation*}

The Lemma is now proved. $\Box$
\vskip 5pt\noindent {\bf Proof of Lemma \ref{le2.2}.} By Lemmas \ref{le2.1} and \ref{le1}, we get that
\begin{equation*}
\begin{split}
\displaystyle \bigg|\int^1_0\mathcal{T}_{\alpha,\beta}(t,x,y)v_1(y)dy\bigg|_{p}^2
&= \sum\limits_{k=1}^{\infty}\lambda_k^{p}
|E_{\alpha, 1}(-\lambda_k^\beta t^\alpha)|^2(v_1, e_k)^2\\
 &\leq t^{\frac{\alpha(q-p)}{\beta}}\sum\limits_{k=1}^{\infty}
\frac{C(\lambda_k^{\beta}t^{\alpha})^{\frac{p-q}{\beta}}}{(1+\lambda_k^\beta t^\alpha)^2}
\lambda_k^{q}(v_1, e_k)^2 \leq C t^{\frac{\alpha(q-p)}{\beta}}|v_1|_{q}^2,
\end{split}
\end{equation*}
where we have used $\frac{(\lambda_k^{\beta}t^{\alpha})
^{\frac{p-q}{\beta}}}{(1+\lambda_k^\beta t^\alpha)^2} \leq C$ for $0\leq q \leq p \leq  2\beta$.

Similarly, we deduce that
\begin{equation*}
\begin{split}
\displaystyle\bigg|\int^1_0\mathcal{R}_{\alpha,\beta}(t,x,y)v_2(y)dy\bigg|_{p}^2
\leq C t^{2-\frac{\alpha(p-r)}{\beta}}|v_2|_{r}^2.
\end{split}
\end{equation*}

Now we consider the case $q>p$. Since $0<\lambda_1 \leq \lambda_2 \leq \lambda_3\leq \cdots$ and $\lambda_k \rightarrow \infty$ as $k\rightarrow \infty$, we obtain from Lemmas \ref{le2.1} and \ref{le1} that
\begin{equation*}
\begin{split}
\displaystyle \bigg|\int^1_0\mathcal{T}_{\alpha,\beta}(t,x,y)v_1(y)dy\bigg|_{p}^2 & \leq \sum\limits_{k=1}^{\infty}
\frac{C}{\lambda_k^{q-p}(1+\lambda_k^\beta t^\alpha)^2}
\lambda_k^{q}(v_1, e_k)^2 \leq Ct^{-2\alpha} |v_1|_{q}^2,
\end{split}
\end{equation*}
and in the similar way, we have
\begin{equation*}
\begin{split}
\displaystyle \bigg|\int^1_0\mathcal{R}_{\alpha,\beta}(t,x,y)v_2(y)dy\bigg|_{p}^2
\leq C t^{2-2\alpha}|v_2|_{r}^2.
\end{split}
\end{equation*}
Thus, \eqref{eq2.11} follows immediately by the triangle inequality.

On the other hand, it follows from Lemmas \ref{le2.1} and \ref{le1} that
\begin{equation*}
\begin{split}
\displaystyle\bigg|\partial_t^\alpha \int_0^1\mathcal{T}_{\alpha,\beta}(t,x,y)v_1(y)dy\bigg|_{p}^2
&= \sum\limits_{l=1}^{\infty}\lambda_l^{p+2\beta}\left(\int_0^1\mathcal{T}_{\alpha,\beta}(t,x,y)v_1(y)dy, e_l\right)^2\\
\displaystyle & = \sum\limits_{l=1}^{\infty}\lambda_l^{p+2\beta}
|E_{\alpha, 1}(-\lambda_l^\beta t^\alpha)|^2(v_1, e_l)^2\\
\displaystyle & \leq t^{-\alpha(2+\frac{p-q}{\beta})}\sum\limits_{l=1}^{\infty}
\frac{C(\lambda_l^{\beta}t^{\alpha})^{\frac{2\beta+p-q}{\beta}}}{(1+\lambda_l^\beta t^\alpha)^2}
\lambda_l^{q}(v_1, e_l)^2\\
\displaystyle & \leq C t^{-\alpha(2+\frac{p-q}{\beta})}|v_1|_{q}^2.
\end{split}
\end{equation*}
A similar estimate for $|\partial_t^\alpha \int_0^1\mathcal{R}_{\alpha,\beta}(t,x,y)v_2(y)dy|_{p}^2$ holds, and this completes the proof. $\Box$
\vskip 5pt\noindent {\bf Proof of Theorem \ref{thm2.1}.} Without loss of generality, we assume that there exists a positive integer ${I_t}$ such that $t=t_{I_t+1}$.
Subtracting \eqref{eq2.8} from \eqref{eq2.7}, we have
\begin{equation}\label{A1}
\begin{split}
\displaystyle u(t,x)-u_n(t,x) &= \int^t_0\int^1_0(t-s)^{\alpha-1}\sum^{\infty}_{k=1}
{E}_{\alpha,\alpha}(-\lambda^{\beta}_k(t-s)^{\alpha})e_k(y)e_k(x)dW(s,y)\\
&\relphantom{=}{}-\int^t_0\int^1_0(t-s)^{\alpha-1}\sum^{\infty}_{k=1}
{E}_{\alpha,\alpha}(-\lambda^{\beta}_k(t-s)^{\alpha})e_k(y)e_k(x)dW_n(s,y),
\end{split}
\end{equation}
where
\begin{equation}\label{w1}
\displaystyle dW(s,y)=\frac{\partial^2W}{\partial s\partial y}dyds
=\sum^{\infty}_{k=1}\sigma_k(s)e_k(y)dyd{\xi}_k(s),
\end{equation}
and
\begin{equation}\label{w3}
\displaystyle dW_n(s,y)=\frac{\partial^2{W_n}}{\partial s\partial y}dyds
=\sum^{\infty}_{k=1}\sigma^n_k(s)e_k(y)\left(\sum\limits_{i=1}^{I_t}
\frac{1}{\sqrt{\Delta t}}\xi_{ki}\chi_i(s)\right)dyds.
\end{equation}
To estimate \eqref{A1}, we introduce an intermediate noise form
\begin{equation}\label{w2}
\displaystyle d\hat{W}_n(s,y)=\frac{\partial^2\hat{W}_n}{\partial s\partial y}dyds
=\sum^{\infty}_{k=1}\sigma^n_k(s)e_k(y)dyd{\xi}_k(s).
\end{equation}
Let
\begin{equation}\label{A13}
\begin{split}
\displaystyle \mathcal{G}_1
\displaystyle  &= \int^t_0\int^1_0(t-s)^{\alpha-1}\sum^{\infty}_{k=1}
{E}_{\alpha,\alpha}(-\lambda^{\beta}_k(t-s)^{\alpha})e_k(y)e_k(x)dW(s,y)\\
&\relphantom{=}{}-\int^t_0\int^1_0(t-s)^{\alpha-1}\sum^{\infty}_{k=1}
{E}_{\alpha,\alpha}(-\lambda^{\beta}_k(t-s)^{\alpha})e_k(y)e_k(x)d\hat{W}_n(s,y),
\end{split}
\end{equation}
and
\begin{equation}\label{A12}
\begin{split}
\displaystyle \mathcal{G}_2&=\int^t_0\int^1_0(t-s)^{\alpha-1}\sum^{\infty}_{k=1}{E}_{\alpha,\alpha}
(-\lambda^{\beta}_k(t-s)^{\alpha})e_k(y)e_k(x)d\hat{W}_n(s,y)\\
&\relphantom{=}{}-\int^t_0\int^1_0(t-s)^{\alpha-1}\sum^{\infty}_{k=1}
{E}_{\alpha,\alpha}(-\lambda^{\beta}_k(t-s)^{\alpha})e_k(y)e_k(x)dW_n(s,y).
\end{split}
\end{equation}
Then it follows from \eqref{A1} and \eqref{A13}-\eqref{A12} that
\begin{equation}\label{A10}
\begin{split}
\displaystyle u-u_n=\mathcal{G}_1+\mathcal{G}_2.
\end{split}
\end{equation}
For $\mathcal{G}_1$, since $\{e_k\}_{k=1}^\infty$ is an orthonormal basis in $L^2(0,1)$, we get from \eqref{w1},  \eqref{w2} and the It\^{o} isometry that
\begin{equation}\label{A5}
\begin{split}
\mathbf{E}\|\mathcal{G}_1\|^2&
= \mathbf{E}\int^1_0\Bigg(\int^t_0\int^1_0(t-s)^{\alpha-1}\sum^{\infty}_{k=1}
{E}_{\alpha,\alpha}(-\lambda^{\beta}_k(t-s)^{\alpha})e_k(y)e_k(x)dW(s,y)\\
&\relphantom{=}{}
-\int^t_0\int^1_0(t-s)^{\alpha-1}\sum^{\infty}_{k=1}
{E}_{\alpha,\alpha}(-\lambda^{\beta}_k(t-s)^{\alpha})e_k(y)e_k(x)d\hat{W}_n(s,y)\Bigg)^2dx\\
\displaystyle & = \int^1_0\mathbf{E}\bigg(\int^t_0(t-s)^{\alpha-1}\sum^{\infty}_{k=1}
{E}_{\alpha,\alpha}(-\lambda^{\beta}_k(t-s)^{\alpha})e_k(x)(\sigma_k(s)-\sigma^n_k(s))d\xi_k(s)\bigg)^2dx\\
\displaystyle & = \int^1_0\int^t_0(t-s)^{2\alpha-2}\left(\sum^{\infty}_{k=1}
{E}_{\alpha,\alpha}(-\lambda^{\beta}_k(t-s)^{\alpha})e_k(x)(\sigma_k(s)-\sigma^n_k(s))\right)^2dsdx.
\end{split}
\end{equation}
Using Lemma \ref{le2.1}, in view of $1<\alpha<2$, we obtain
\begin{equation}\label{A6}
\begin{split}
\mathbf{E}\|\mathcal{G}_1\|^2
\displaystyle & \leq \sum^{\infty}_{k=1}(\eta_k^n)^2\int^t_0(t-s)^{2\alpha-2}
{E}^2_{\alpha,\alpha}(-\lambda^{\beta}_k(t-s)^{\alpha})ds\\
\displaystyle & \leq C\sum^{\infty}_{k=1}\lambda_k^{-\frac{2\beta(\alpha-1)}{\alpha}}(\eta_k^n)^2\int^t_0
\bigg|\frac{(\lambda^{\beta}_k(t-s)^{\alpha})^{\frac{\alpha-1}{\alpha}}}
{1+\lambda^{\beta}_k(t-s)^{\alpha}}\bigg|^2ds\\
\displaystyle & \leq C\sum^{\infty}_{k=1}\lambda_k^{-\frac{2\beta(\alpha-1)}{\alpha}}(\eta_k^n)^2,
\end{split}
\end{equation}
where we have used $\bigg|\frac{(\lambda^{\beta}_k(t-s)^{\alpha})^{\frac{\alpha-1}{\alpha}}}
{1+\lambda^{\beta}_k(t-s)^{\alpha}}\bigg|^2\leq C$.

For  $\mathcal{G}_2$, since  $\{e_k\}_{k=1}^\infty$ is an orthonormal basis in $L^2(0,1)$, then by \eqref{w3}-\eqref{w2} and \eqref{A12}, we have
\begin{equation}\label{A7}
\begin{split}
\displaystyle \mathbf{E}\|\mathcal{G}_2\|^2 &= \mathbf{E}
\int_0^1\Bigg(\int_0^t(t-s)^{\alpha-1}\sum\limits_{k=1}^{\infty}\int_0^1
E_{\alpha,\alpha}(-\lambda_k^\beta (t-s)^\alpha)e_k(y)e_k(x)d\hat{W}_n(s,y) \\
&\relphantom{=}{}-\int_0^t(t-s)^{\alpha-1}\sum\limits_{k=1}^{\infty}\int_0^1
E_{\alpha,\alpha}(-\lambda_k^\beta (t-s)^\alpha)e_k(y)e_k(x)d{W}_n(s,y)\Bigg)^2dx\\
&=\mathbf{E}\int_0^1\Bigg(\int_0^t(t-s)^{\alpha-1}\sum\limits_{k=1}^{\infty}
E_{\alpha,\alpha}(-\lambda_k^\beta (t-s)^\alpha)e_k(x)\sigma_k^n(s)\\
&\relphantom{=============}{}\times\left(\dot{\xi}_k(s)-\sum\limits_{i=1}^{I_t}
\frac{1}{\sqrt{\Delta t}}\xi_{ki}\chi_i(s)\right)ds\Bigg)^2dx\\
&=\int_0^1\mathbf{E}\Bigg(\sum\limits_{l=1}^{I_t}
\int_{t_{l}}^{t_{l+1}}\bigg((t-s)^{\alpha-1}\sum\limits_{k=1}^{\infty}
E_{\alpha,\alpha}(-\lambda_k^\beta (t-s)^\alpha)e_k(x)\sigma_k^n(s)\\
&\relphantom{=}{}-\frac{1}{\Delta t}\int_{t_{l}}^{t_{l+1}}(t-\tilde{s})^{\alpha-1}\sum\limits_{k=1}^{\infty}
E_{\alpha,\alpha}(-\lambda_k^\beta (t-\tilde{s})^\alpha)
e_k(x)\sigma_k^n(\tilde{s})d\tilde{s}\bigg)d\xi_k(s)\Bigg)^2dx.
\end{split}\end{equation}
Using the It\^{o} isometry, we get
\begin{equation}\label{A11}
\begin{split}
\displaystyle\mathbf{E}\|\mathcal{G}_2\|^2 &=\int_0^1\sum\limits_{l=1}^{I_t}
\int_{t_{l}}^{t_{l+1}}\bigg((t-s)^{\alpha-1}\sum\limits_{k=1}^{\infty}
E_{\alpha,\alpha}(-\lambda_k^\beta (t-s)^\alpha)e_k(x)\sigma_k^n(s)\\
&\relphantom{=}{}-\frac{1}{\Delta t}\int_{t_{l}}^{t_{l+1}}(t-\tilde{s})^{\alpha-1}\sum\limits_{k=1}^{\infty}
E_{\alpha,\alpha}(-\lambda_k^\beta (t-\tilde{s})^\alpha)
e_k(x)\sigma_k^n(\tilde{s})d\tilde{s}\bigg)^2dsdx\\
&=\int_0^1\sum\limits_{l=1}^{I_t}\frac{1}{(\Delta t)^2}
\int_{t_{l}}^{t_{l+1}}\Bigg(\int_{t_{l}}^{t_{l+1}}(t-s)^{\alpha-1}\sum\limits_{k=1}^{\infty}
E_{\alpha,\alpha}(-\lambda_k^\beta (t-s)^\alpha)e_k(x)\\
&\relphantom{==}{}\times(\sigma_k^n(s)-\sigma_k^n(\tilde{s}))d\tilde{s}+\int_{t_{l}}^{t_{l+1}}
\sum\limits_{k=1}^{\infty}\bigg((t-{s})^{\alpha-1}
E_{\alpha,\alpha}(-\lambda_k^\beta (t-{s})^\alpha)\\
&\relphantom{=======}{}-(t-\tilde{s})^{\alpha-1}E_{\alpha,\alpha}(-\lambda_k^\beta (t-\tilde{s})^\alpha)\bigg)
e_k(x)\sigma_k^n(\tilde{s})d\tilde{s}\Bigg)^2dsdx.
\end{split}\end{equation}
By Lagrange's mean value theorem, we have
\begin{equation}\label{A8}
\begin{split}
\displaystyle\mathbf{E}\|\mathcal{G}_2\|^2 & =\int_0^1 \sum\limits_{l=1}^{I_t}\frac{1}{(\Delta t)^2}
\int_{t_{l}}^{t_{l+1}}\bigg(\int_{t_{l}}^{t_{l+1}}(t-s)^{\alpha-1}\sum\limits_{k=1}^{\infty}
E_{\alpha,\alpha}(-\lambda_k^\beta (t-s)^\alpha)\\
&\relphantom{=====================}{}\times e_k(x)(\sigma_k^{n})^\prime(\zeta_l)(s-\tilde{s})d\tilde{s}\bigg)^2dsdx\\
&\relphantom{=}{}+\int_0^1\sum\limits_{l=1}^{I_t}\frac{1}{(\Delta t)^2}
\int_{t_{l}}^{t_{l+1}}\bigg(\int_{t_{l}}^{t_{l+1}}\sum\limits_{k=1}^{\infty}\bigg((t-{s})^{\alpha-1}
E_{\alpha,\alpha}(-\lambda_k^\beta (t-{s})^\alpha)\\
&\relphantom{=========}{}-(t-\tilde{s})^{\alpha-1}
E_{\alpha,\alpha}(-\lambda_k^\beta (t-\tilde{s})^\alpha)\bigg)e_k(x)\sigma_k^n(\tilde{s})d\tilde{s}\bigg)^2dsdx\\
&:=\mathcal{I}_1+\mathcal{I}_2,
\end{split}\end{equation}
where $\zeta_l$ is between $s$ and $\tilde{s}$.

Note that $\{e_k\}_{k=1}^\infty$ is an orthonormal basis in $L^2(0,1)$.
For $\mathcal{I}_1$, using the smoothness assumption on $\sigma^{n}_k(t)$, then it follows from Lemma \ref{le2.1} that
\begin{equation}\label{A9}
\begin{split}
\mathcal{I}_1&\leq \int_0^1\sum\limits_{l=1}^{I_t}\frac{1}{(\Delta t)^2}
\int_{t_{l}}^{t_{l+1}}\left((t-s)^{\alpha-1}\sum\limits_{k=1}^{\infty}E_{\alpha,\alpha}(-\lambda_k^\beta (t-s)^\alpha)e_k(x)\gamma_k^n(\Delta t)^2\right)^2 dsdx\\
&\leq (\Delta t)^2\sum\limits_{l=1}^{I_t}\int_{t_{l}}^{t_{l+1}}(t-s)^{2\alpha-2}
\sum\limits_{k=1}^{\infty}E^2_{\alpha,\alpha}(-\lambda_k^\beta (t-s)^\alpha)(\gamma_k^n)^2ds\\
&\leq C(\Delta t)^2\int^t_0\sum\limits_{k=1}^{\infty}\lambda_k^{-\frac{2\beta(\alpha-1)}{\alpha}}(\gamma_k^n)^2
\bigg|\frac{(\lambda^{\beta}_k(t-s)^{\alpha})^{\frac{\alpha-1}{\alpha}}}
{1+\lambda^{\beta}_k(t-s)^{\alpha}}\bigg|^2ds\\
&\leq C(\Delta t)^2\sum\limits_{k=1}^{\infty}\lambda_k^{-\frac{2\beta(\alpha-1)}{\alpha}}(\gamma_k^n)^2,
\end{split}
\end{equation}
where we have used $\frac{(\lambda^{\beta}_k(t-s)^{\alpha})^{\frac{\alpha-1}{\alpha}}}
{1+\lambda^{\beta}_k(t-s)^{\alpha}}\leq C$.

For $\mathcal{I}_2$, we need the following equality (see \cite[(1.82)]{Podlubny})
\[\frac{d}{d\tau}(t-\tau)^{\alpha-1}E_{\alpha,\alpha}(-\lambda^{\beta}_k(t-\tau)^{\alpha})
=-(t-\tau)^{\alpha-2}E_{\alpha,\alpha-1}(-\lambda^\beta_k(t-\tau)^\alpha).\]
By Lemma \ref{le2.1}, we deduce that
\begin{equation*}
\begin{split}
\displaystyle &(t-{s})^{\alpha-1}E_{\alpha,\alpha}(-\lambda_k^\beta (t-{s})^\alpha)-(t-\tilde{s})^{\alpha-1}
E_{\alpha,\alpha}(-\lambda_k^\beta (t-\tilde{s})^\alpha)\\
&\relphantom{=}{}=\int_{\tilde{s}}^{s}(t-\tau)^{\alpha-2}(-E_{\alpha,\alpha-1}(-\lambda^\beta_k(t-\tau)^\alpha))d\tau
\leq C\int_{\tilde{s}}^{s}(t-\tau)^{\alpha-2}d\tau:=\mathcal{I}_2^1
\end{split}
\end{equation*}
Note that $a^\theta-b^\theta\leq (a-b)^\theta$ for $a>b>0$ and $0<\theta<1$. Let $0<\varepsilon<\frac{1}{2}$ be given arbitrary. For $\displaystyle 1<\alpha\leq \frac{3}{2}$, we have
\begin{equation}\label{A2}
\begin{split}\displaystyle
\mathcal{I}_2^1&=C\int_{\tilde{s}}^{s}(t-\tau)^{-\frac{1}{2}+\varepsilon}(t-\tau)^{\alpha-\frac{3}{2}-\varepsilon}d\tau\\
&\leq C(t-s)^{-\frac{1}{2}+\varepsilon}\int_{\tilde{s}}^{s}(t-\tau)^{\alpha-\frac{3}{2}-\varepsilon}d\tau
\leq C(t-s)^{-\frac{1}{2}+\varepsilon}(\Delta t)^{\alpha-\frac{1}{2}-\varepsilon}.
\end{split}
\end{equation}
Note that $\{e_k\}_{k=1}^\infty$ is an orthonormal basis in $L^2(0,1)$. For $\mathcal{I}_2$, by the assumption on $\sigma^{n}_k(t)$ and \eqref{A2}, we get
\begin{equation}\label{A3}
\begin{split}\displaystyle
\mathcal{I}_2&\leq C\int_0^1 \sum\limits_{l=1}^{I_t}
\frac{1}{(\Delta t)^2}\int_{t_{l}}^{t_{l+1}}\bigg(\int_{t_{l}}^{t_{l+1}}\sum^{\infty}_{k=1}
(t-s)^{-\frac{1}{2}+\varepsilon}(\Delta t)^{\alpha-\frac{1}{2}-\varepsilon}e_k(x)\mu_k^nd\tilde{s}\bigg)^2dsdx\\
&\leq C\int_0^1\int_{0}^{t}(t-s)^{-1+2\varepsilon}(\Delta
t)^{2\alpha-1-2\varepsilon}\left(\sum^{\infty}_{k=1}e_k(x)\mu_k^n\right)^2dsdx\\
&\leq C(\Delta t)^{2\alpha-1-2\varepsilon}\sum^{\infty}_{k=1}(\mu_k^n)^2 \int_0^t (t-s)^{-1+2\varepsilon}ds\\
&\leq C t^{2\varepsilon} (\Delta t)^{2\alpha-1-2\varepsilon}\sum^{\infty}_{k=1}(\mu_k^n)^2.
\end{split}
\end{equation}
For $\frac{3}{2}<\alpha<2$, we have $\mathcal{I}_2^1\leq C(t-s)^{\alpha-2}\Delta t$. Arguing as in the argument of \eqref{A3}, we deduce that
\begin{equation}\label{A4}
\begin{split}
\mathcal{I}_2&\leq C\int_0^1\sum\limits_{l=1}^{I_t}
\frac{1}{(\Delta t)^2}\int_{t_{l}}^{t_{l+1}}\bigg(\int_{t_{l}}^{t_{l+1}}\sum^{\infty}_{k=1}
(t-s)^{\alpha-2}\Delta t e_k(x)\mu_k^nd\tilde{s}\bigg)^2dsdx\\
&\leq C(\Delta t)^2 \int_0^1\int_0^t (t-s)^{2\alpha-4}\left(\sum^{\infty}_{k=1}e_k(x)\mu_k^n \right)^2dsdx\\
&\leq C(\Delta t)^2 \sum^{\infty}_{k=1}(\mu_k^n)^2\int_0^t (t-s)^{2\alpha-4}ds\leq Ct^{2\alpha-3}(\Delta t)^2 \sum^{\infty}_{k=1}(\mu_k^n)^2.
\end{split}
\end{equation}
For $1<\alpha\leq\frac{3}{2}$, by \eqref{A10}, \eqref{A6}, \eqref{A8}, \eqref{A9} and \eqref{A3}, we conclude that for
any $0<\varepsilon<\frac{1}{2}$,
\begin{equation*}
\begin{split}
\mathbf{E}\|u(t)-u_n(t)\|^2&\leq C\sum^{\infty}_{k=1}\lambda_k^{-\frac{2\beta(\alpha-1)}{\alpha}}(\eta_k^n)^2
+C(\Delta t)^2\sum^{\infty}_{k=1}\lambda_k^{-\frac{2\beta(\alpha-1)}{\alpha}}(\gamma_k^n)^2\\
&\relphantom{=}{}
+Ct^{2\varepsilon}(\Delta t)^{2\alpha-1-2\varepsilon}\sum^{\infty}_{k=1}(\mu_k^n)^2, \quad t>0;
\end{split}
\end{equation*}
and for $\frac{3}{2}< \alpha<2$, we obtain from \eqref{A10}, \eqref{A6}, \eqref{A8}, \eqref{A9} and \eqref{A4} that
\begin{equation*}
\begin{split}
\mathbf{E}\|u(t)-u_n(t)\|^2&\leq C\sum^{\infty}_{k=1}\lambda_k^{-\frac{2\beta(\alpha-1)}{\alpha}}(\eta_k^n)^2
+C(\Delta t)^2\sum^{\infty}_{k=1}\lambda_k^{-\frac{2\beta(\alpha-1)}{\alpha}}(\gamma_k^n)^2\\
&\relphantom{=}{}+Ct^{2\alpha-3}(\Delta t)^2\sum^{\infty}_{k=1}(\mu_k^n)^2,\quad t>0.
\end{split}
\end{equation*}
The proof is completed. $\Box$
\vskip 5pt\noindent {\bf Proof of Lemma \ref{le3.2}.}
By Lemmas \ref{le2.1} and \ref{le1}, we obtain that
\begin{equation}\label{6.5}
\begin{split}
\displaystyle &\bigg|\int^1_0\mathcal{S}^h_{\alpha,\beta}(t,x,y)\psi(y)dy\bigg|_{p,h}^2
= \sum\limits_{k=1}^{N}(\lambda_k^{h,\beta})^{\frac{p}{\beta}} t^{2\alpha-2}
|E_{\alpha, \alpha}(-\lambda_k^{h,\beta} t^\alpha)|^2(\psi, e_k^h)^2\\
\displaystyle &\relphantom{=}{}  \leq t^{\frac{\alpha(q-p)}{\beta}-2+2\alpha}\max\limits_{k}
\frac{C\left(\lambda_k^{h,\beta}t^{\alpha}\right)^{\frac{p-q}{\beta}}}{(1+\lambda_k^{h,\beta} t^\alpha)^2}
\sum\limits_{k=1}^{N}(\lambda_k^{h,\beta})^{\frac{q}{\beta}}(\psi, e_k^h)^2\\
\displaystyle &\relphantom{=}{}  \leq C t^{\frac{\alpha(q-p)}{\beta}-2+2\alpha}|\psi|_{q, h}^2,
\end{split}
\end{equation}
where we have used $\frac{\left(\lambda_k^{h,\beta}t^{\alpha}\right)^{\frac{p-q}{\beta}}}{(1+\lambda_k^{h,\beta} t^\alpha)^2}\leq C$ for $ p - 2\beta \leq q \leq p$.

For $q>p$, since $\{\lambda^{h,\beta}_k\}$ are bounded away from zero independent of the mesh size $h$, we deduce from Lemmas \ref{le2.1} and \ref{le1} that
\begin{equation}\label{6.6}
\begin{split}
\displaystyle \bigg|\int^1_0\mathcal{S}^h_{\alpha,\beta}(t,x,y)\psi(y)dy\bigg|_{p,h}^2
& \leq t^{2\alpha-2}\max\limits_{k}
\frac{C}{(1+\lambda_k^{h,\beta} t^\alpha)^2 (\lambda_k^{h,\beta})^{\frac{q-p}{\beta}}}
\sum\limits_{k=1}^{N}(\lambda_k^{h,\beta})^{\frac{q}{\beta}}(\psi, e_k^h)^2\\
& \leq C t^{-2}|\psi|_{q, h}^2.
\end{split}
\end{equation}
Thus, we complete the proof of Lemma \ref{le3.2}. $\Box$

\bibliographystyle{amsplain}

\end{document}